\documentclass[
11pt,reqno]{amsart}

\usepackage{amsthm, amsfonts, amssymb}

\theoremstyle{definition}
\newtheorem{ntn}{Notation}[section]
\newtheorem{dfn}[ntn]{Definition}
\theoremstyle{plain}
\newtheorem{lem}[ntn]{Lemma}
\newtheorem{prp}[ntn]{Proposition}
\newtheorem{thm}[ntn]{Theorem}
\newtheorem{cor}[ntn]{Corollary}
\newtheorem{conj}[ntn]{Conjecture}
\theoremstyle{remark}
\newtheorem{rem}{Remark}
\newtheorem{exa}{Example}

\def\floor[#1]{\lfloor #1 \rfloor }

\newcommand{\z}{\mathbb{Z}}

\newcommand{\q}{\mathbb{Q}}
\newcommand{\R}{\mathbb{R}}
\newcommand{\lan}{\langle}
\newcommand{\ran}{\rangle}

\newcommand{\GL}{\mathit{{\rm GL}}}
\newcommand{\SL}{\mathit{{\rm SL}}}

\newcommand{\EE}{{\mathcal{E}}}
\newcommand{\DD}{{\mathfrak{d}}}
\newcommand{\ppp}{\mathfrak{p}}

\newcommand{\ee}{\mathcal{E}}

\newcommand{\inc}{{\rm inc}}
\newcommand{\id}{{\rm id}}

\newcommand{\tors}{{{\rm Tor}_1^{\z}}}
\newcommand{\zzz}{\z[\frac{1}{2}]}

\newcommand{\s}{\Sigma}
\newcommand{\si}{\sigma}

\newcommand{\arr}{\rightarrow}
\newcommand{\larr}{\longrightarrow}

\newcommand{\se}{\subseteq}

\newcommand{\mt}{\mapsto}
\newcommand{\two}{\twoheadrightarrow}

\newcommand{\fn}{F^n}

\newcommand{\fff}{{F^\ast}}

\newcommand{\eee}{{\tilde{E}}}

\newcommand{\rr}{{F^\ast}}

\newcommand{\stabe}{{\rm Stab}}
\renewcommand{\char}{{\rm char}}
\newcommand{\diag}{{\rm diag}}
\renewcommand{\ker}{{\rm ker}}
\newcommand{\coker}{{\rm coker}}
\newcommand{\im}{{\rm im}}
\newcommand{\ind}{{\rm ind}}

\newcommand{\nnn}{\z[\frac{1}{(n-1)!}]}
\newcommand {\mtx}[4]
{\left(
\begin{array}{cc}
#1 & #2   \\
#3 & #4
\end{array}
\right)}

\newtheoremstyle{athm}
  {}
  {}
  {\itshape}
  {}
  {\scshape}
  {}
  {.5em}
  {\thmnote{#3}}
\theoremstyle{athm}
\newtheorem*{athm}{}

\begin{document}

\title{Homology of $\SL_n$ and $\GL_n$
over an infinite field}
\author{B. Mirzaii}

\begin{abstract}
The homology of $\GL_n(F)$ and $\SL_n(F)$
is studied, where $F$ is an infinite field.
Our main theorem states that the natural map
$H_4(\GL_3(F), k) \arr H_4(\GL_4(F), k)$
is injective where $k$ is a field with
$\char(k) \neq 2, 3$. For  algebraically closed
field $F$, we prove a better result, namely,
$H_4(\GL_3(F), \z) \arr H_4(\GL_4(F), \z)$ is
injective. We will prove a similar result
replacing $\GL$ by $\SL$. This is used to investigate 
the indecomposable part of the $K$-group
$K_4(F)$.
\end{abstract}

\maketitle

\section{Introduction}

In the beginning of the 1970's two type of $K$-groups in algebra
appeared: Quillen's $K$-groups and Milnor's $K$-groups. For a
field $F$, Quillen defined the $K$-group $K_n(F)$ as the $n$-th
homotopy group of the space $B\GL(F)^+$ and Milnor defined the
$K$-group $K_n^M(F)$ as the $n$-th degree part of $T(\fff)/\lan a
\otimes (1-a):a \in \fff-\{1\}\ran$, where $T(\fff):=\z \oplus
\fff \oplus \fff \otimes\fff \oplus \cdots$ is the tensor algebra
of $\fff$. There is a canonical ring homomorphism $K_\ast^M(F)
\arr K_\ast(F)$, therefore a canonical homomorphism $K_n^M(F) \arr
K_n(F)$. On the other hand, the Hurewicz theorem, in algebraic
topology, relates homotopy groups to homology groups, which are
much easier to calculate. This in turn provides a homomorphism
from $K_n(F)$ to the $n$-th integral homology of the stable group
$\GL(F)$.

One of the important approaches to investigate $K$-groups is by
means of their relation with  integral homology groups of $\GL(F)$
and Milnor $K$-groups.
Suslin's stability theorem states that for an infinite field  $F$
the natural map $H_i(\GL_n(F),\z) \arr H_i(\GL(F),\z)$ is bijective if
$n \ge i$ \cite{suslin1985}. Using this 
Suslin  constructed a map from
$H_n(\GL_n(F), \z)$ to Milnor's $K_n$-group $K_n^M(F)$, denoted by $s_n$, 
such that the sequence
\begin{gather*}\label{suslinexactness}
H_n(\GL_{n-1}(F), \z) 
\overset{H_n(\inc)}{\larr} H_n(\GL_n(F), \z)
\overset{s_n}{\larr} K_n^M(F) \larr 0
\end{gather*}
is exact. Combining these two results he constructed a map from
$K_n(F)$ to  $K_n^M(F)$ such that the composite homomorphism
\[
K_n^M(F) \arr K_n(F) \arr K_n^M(F)
\]
coincides with the multiplication by $(-1)^{n-1}(n-1)!$ 
\cite[Sec. 4]{suslin1985}.

Now one might ask about the kernel of $H_n(\inc)$ in the above
exact sequence. In this direction,
Suslin posed a problem, which is now
referred to as `a conjecture by Suslin' (see \cite[4.13]{sah1989}
and \cite[7.7]{borel-yang1994}).

\begin{athm}[{\bf Injectivity Conjecture.}]\label{conj1}
For any infinite field $F$ the natural homomorphism
\[
H_n(\inc):H_n(\GL_{n-1}(F), \q) \arr H_n(\GL_n(F), \q)
\]
is injective.
\end{athm}

This conjecture is easy if $n=1,2$.
For $n=3$ the conjecture was proved positively by Sah \cite{sah1989} and
Elbaz-Vincent \cite{elbaz1998}.
The conjecture is proven in full for number fields \cite{borel-yang1994}.
The proof of this conjecture for $n=4$ is
the main goal of this paper (Theorem \ref{g4-inj}).

Here we take a general step towards this conjecture.
We show that up to an induction step, the above conjecture follows from the
exactness of a certain complex. Let 
\[
P_n:=H_n(inc): H_n(\GL_{n-1}(F), \q)\arr H_n(\GL_{n}(F), \q)
\]
and let $Q_n$ be the complex
\begin{gather*}\label{exactmirzaii}
\hspace{-3 cm}
H_n(\fff^2\times \GL_{n-2}(F), \q)
\overset{\beta_2^{(n)}}{\arr}
H_n(\fff \times \GL_{n-1}(F), \q)\\
\hspace{7 cm}
\overset{\beta_1^{(n)}}{\arr}
H_n(\GL_n(F), \q) \arr 0,
\end{gather*}
where $\beta_1^{(n)}=H_n(\inc)$ and
$\beta_2^{(n)}=H_n(\alpha)-H_n(\inc)$, \[
\alpha:\fff^2\times
\GL_{n-2}(F) \arr\fff \times \GL_{n-1}(F),\ \ 
\diag(a, b, A) \arr \diag(b, a, A).
\]


\begin{athm}[{\bf Proposition.}] Let $F$ be an infinite field.
If $Q_n$ is exact and if $P_{n-2}$ and $P_{n-1}$ are injective,
then $P_n$ is injective.
%
\end{athm}
The proof follows from a careful analysis of some spectral sequences,
connecting the homology of the groups $\fff^p \times \GL_{n-p}(F)$ 
for different values of $p$.
One of the main ingredient in the proof of this proposition
is a construction of an explicit map from
$K_n^M(F)$ to $H_n(\GL_n(F), \z)$, denoted by $\nu_n$
(this was only known for $n=2$).
This construction fits in our previous theory, that is,  the
composite homomorphism
\[
K_n^M(F) \overset{\nu_n}{\larr} 
H_n(\GL_n(F), \z)
\overset{s_n}{\larr} K_n^M(F)
\]
coincides with the multiplication by $(-1)^{n-1}(n-1)!$.
As we mentioned above, there is a homomorphism from
$K_n^M(F)$ to $H_n(\GL_n(F), \z)$ that factors through $K_n(F)$,
denoted by $h_n$. We don't know whether these two maps coincide.

Here are our main results.

\begin{athm}[{\bf Theorem.}] Let F be an infinite field.
\par {\rm (i)} The complex
\begin{gather*}
H_4(\fff^2\times \GL_{2}(F), \z)
\overset{\beta_2^{(4)}}{\arr}
H_4(\fff \times \GL_{3}(F), \z)
\overset{\beta_1^{(4)}}{\arr}
H_4(\GL_4(F), \z) \arr 0
\end{gather*}
is exact.
\par {\rm (ii)} Let $k$ be a field with $\char(k) \neq 2, 3$. Then
\[
H_4(\inc):H_4(\GL_{3}(F), k) \arr H_4(\GL_4(F), k)
\]
is injective.
\par {\rm (iii)} If $F$ is algebraically closed, then
\[
H_4(\inc):H_4(\GL_{3}(F), \z) \arr H_4(\GL_4(F), \z)
\] is injective.
\end{athm}


Also we show that a similar result as in part (ii) and (iii)
of the above theorem is true if we replace $\GL$ with $\SL$.
As an application we will study the indecomposable part
of $K_4(F)$. Namely we will prove that for an algebraically
closed field $F$, $K_4(F)^\ind:=\coker(K_4^M(F)
{\arr} K_4(F))$ embeds in $H_4(\SL_3(F), \z)$.


Here we establish some notation.
In this note, by $H_i(G)$ we mean the $i$-th integral homology of the
group $G$. We use the bar resolution to define the homology of a
group \cite[Chap. I, Section 5]{brown1994}. Define
${\rm \bf{c}}({g}_1, {g}_2,\dots, {g}_n)=
\sum_{\si \in \s_n}
{{\rm sign}(\si)}
[{g}_{\si(1)}| {g}_{\si(2)}|\dots | {g}_{\si(n)}] \in H_n(G)$,
where ${g}_i \in G$ pairwise commute and
$\s_n$ is the symmetric group of degree $n$.
By $\GL_n$ and $\SL_n$  we mean $\GL_n(F)$ and $\SL_n(F)$, where
$F$ is an infinite field. Note that $\GL_0$ is the trivial group and
$\GL_1=\fff$. By $\fff^m$ we mean
$\fff \times \cdots\times \fff$ ($m$-times) or
the subgroup of $\fff$, $\{a^m| a \in \fff\}$, depending on the context.
This shall not cause any confusion.
The $i$-th factor of $\fff^m=\fff \times \cdots\times \fff$, ($m$-times),
is denoted by $F_i^\ast$.


\section{Homology of $\SL_n$}

The action of $\fff$ on $\SL_n$ defined by 
$a. A:=\mtx {a} {0} {0} {1}  A \mtx {a^{-1}} {0} {0} {1}$ induces 
an action of $\fff$ on $H_i(\SL_n)$. In 
this article by $H_i(\SL_n)_\fff$ we simply mean $H_0(\fff, H_i(\SL_n))$. 
It is easy to see that the natural map $\SL_n \arr \SL$ induces 
a map of homology groups $H_i(\SL_n)_\fff \arr H_i(\SL)$.

The map $\delta: \GL \arr \SL$ given by 
$ M \mt \mtx {\det(M)^{-1}} {0} {0} {M}$
induces a homomorphism $H_i(\delta):H_i(\GL) \arr H_i(\SL)$ such that
the composition $H_i(\SL) \overset{H_i(\inc)}{\larr} 
H_i(\GL) \overset{H_i(\delta)}{\larr} H_i(\SL)$ is the identity map.
Therefore $H_i(\SL)$ embeds in  $H_i(\GL)$.


\begin{lem}\label{isom-sl}
{\rm (i)} $ H_i(\SL_n)_\fff \simeq H_i(\SL)$ for  $n \ge i$. 
In particular if $F$ is algebraically closed, then 
$H_i(\SL_n) \simeq H_i(\SL)$ for $n \ge i$.
\par {\rm (ii)} Let $A:=\z[\frac{1}{n-1}]$ and let $q \le n-1$. 
Then for $p \ge 0$ the map
\[
H_p(\fff, H_q(\SL_{n-1}, A)) \arr H_p(\fff, H_q(\SL, A)),
\]
induced by the map of pair
$(\id, \inc): (\fff, \SL_{n-1}) \arr (\fff, \SL)$, is isomorphism.
\end{lem}
\begin{proof}
The part (i) is rather well known (see \cite[2.7]{sah1989}). \\
(ii) If $q=0$, then the claim is trivial. So let $q \ge 1$.
The short exact sequence 
\[
1 \arr \mu_{n-1, F} \arr \fff \overset{(.)^{n-1}}{\arr} 
\fff^{n-1} \arr 1
\]
gives us the
Lyndon-Hochschild-Serre spectral sequence
\begin{gather*}
\ee_{r, s}^2=H_r(\fff^{n-1}, H_s(\mu_{{n-1}, F}, T)) 
\Rightarrow  H_{r+s}(\fff, T),
\end{gather*}
where $T=H_q(\SL_{n-1}, A)$. Since the order of 
$\mu_{{n-1}, F}$ is invertible in $A$,
$H_s(\mu_{{n-1}, F}, T)=0$ for $s \ge 1$ \cite[10.1]{brown1994}. 
Thus $\ee_{r, s}^2=0$ for $s \ge 1$.
The action of $\fff^{n-1}$ on $\mu_{{n-1}, F}$ and $T$ is trivial, so
$\ee_{0, 0}^2= H_0(\mu_{{n-1}, F}, T)=
H_0(\mu_{{n-1}, F}, H_q(\SL_{n-1}, A))$.
From this and (i) one deduces that
\begin{gather*}
\ee_{0, 0}^\infty \simeq H_0(\mu_{{n-1}, F}, H_q(\SL_{n-1}, A))
\simeq H_0(\fff, H_q(\SL_{n-1}, A))\simeq H_q(\SL, A)
\end{gather*}
and therefore
$
\ee_{r, 0}^\infty \simeq \ee_{r, 0}^2
=H_r(\fff^{n-1}, H_q(\SL, A)).
$
An easy analysis shows that
\[
H_r(\fff^{n-1}, H_q(\SL, A)) \simeq H_r(\fff, H_q(\SL_{n-1}, A)).
\]
Once more from the short exact sequence 
$1 \arr \fff^{n-1}  \arr \fff {\arr} \fff / \fff^{n-1} \arr 1$
one gets the Lyndon-Hochschild-Serre spectral sequence
\begin{gather*}
{\ee'}_{r, s}^2=H_r(\fff/ \fff^{n-1}, H_s(\fff^{n-1}, S))
\Rightarrow  H_{r+s}(\fff, S)
\end{gather*}
where $S=H_q(\SL, A)$. It is easy to see that
${\ee'}_{r, s}^2=0$ for $r \ge 1$ and
\[
{\ee'}_{0, s}^\infty \simeq{\ee'}_{0, s}^2=
H_0(\fff/ \fff^{n-1}, H_s(\fff^{n-1}, S))=H_s(\fff^{n-1}, H_q(\SL, A)). 
\]
This implies that
$H_s(\fff^{n-1}, H_q(\SL, A)) \simeq H_s(\fff, H_q(\SL, A))$.
Hence for $r \ge 0$ and $q \le n-1$,
\[
H_r(\fff, H_q(\SL_{n-1}, A)) \simeq H_r(\fff, H_q(\SL, A)).
\]
It is not difficult to see that this isomorphism is 
induced by the map of pair
$({()}^{n-1}, \inc): (\fff, \SL_{n-1}) \arr (\fff, \SL)$. 
In a similar way, one can prove that the map of pair
$({()}^{n-1}, \inc): (\fff, \SL) \arr (\fff, \SL)$
induces the isomorphism 
\[
H_r(\fff, H_q(\SL, A)) \simeq H_r(\fff, H_q(\SL, A)).
\]
Applying the functor $H_r$ to the commutative diagram  
\[
\begin{array}{ccc}
(\fff, \SL_{n-1}) &\overset{(\id, \inc)}
{\!\!\!-\!\!\!-\!\!\!-\!\!\!-\!\!\!-\!\!\!\larr} & (\fff, \SL) \\
\Big\downarrow\vcenter{%
\rlap{$\scriptstyle{(\id, \inc)}$}}
&       & \Big\downarrow\vcenter{%
\rlap{$\scriptstyle{({()}^{n-1}, \inc)}$}} \\
(\fff, \SL_{n-1})  & \overset{({()}^{n-1}, \inc)}
{\!\!\!-\!\!\!-\!\!\!-\!\!\!-\!\!\!-\!\!\!\larr} & (\fff, \SL)
\end{array}
\]
gives us the isomorphism that we are looking for.
\end{proof}

\begin{lem}\label{inter}
{\rm (i)}  If $A:=\z[\frac{1}{n-1}]$, then
\begin{gather*}
\hspace{-2 cm}
\im(H_n(\GL_{n-1}, A) \arr H_n(\GL_n, A)) \cap H_n(\SL_n, A)_\fff \\
\hspace{5 cm}
= \im(H_n(\SL_{n-1}, A)_\fff \arr H_n(\SL_n, A)_\fff).
\end{gather*}
\par {\rm (ii)} If $F$ is algebraically closed, then
\[
\im(H_n(\GL_{n-1}) \arr H_n(\GL_n)) \cap H_n(\SL_n) =
\im(H_n(\SL_{n-1}) \arr H_n(\SL_n)).
\]
\end{lem}
\begin{proof}
To prove (i) consider the associated Lyndon-Hochschild-Serre
spectral sequence of the diagram of extensions
\[
\begin{array}{ccccccccc}
\!\!\! 1  & \!\!\!\larr &\!\!\!  \SL_{n-1}  &
\!\!\!
{\larr}
& \!\!\!  \GL_{n-1} &  \!\!\! \larr
&\!\!\! \fff
& \!\!\!\larr & \!\!\! 1 \\
&  &\!\!\!\Big\downarrow\vcenter{%
\rlap{$\scriptstyle{}$}}
&       & \!\!\!\Big\downarrow\vcenter{%
\rlap{$\scriptstyle{}$}}
&  &\!\!\!\Big\downarrow\vcenter{%
\rlap{$\scriptstyle{}$}}
&      &  \\
\!\!\! 1  & \!\!\!\larr & \!\!\! \SL
& \!\!\!\larr & \!\!\! \GL
& \!\!\!  \larr   & \!\!\! \fff    &\!\!\!
\larr & \!\!\! 1.
\end{array}
\]
{}From this we obtain a map of spectral sequences
\begin{gather*}
\begin{array}{ccc}
E_{p, q}^2=H_p(\fff, H_q(\SL_{n-1}, A)) & 
\Rightarrow & H_{p+q}(\GL_{n-1}, A)\\
\Big\downarrow\vcenter{%
\rlap{$\scriptstyle{}$}} &      &
\Big\downarrow\vcenter{%
\rlap{$\scriptstyle{}$}}       \\
{E'}_{p, q}^2=H_p(\fff, H_q(\SL, A))& \Rightarrow & H_{p+q}(\GL, A).
\end{array}
\end{gather*}
By Lemma \ref{isom-sl} for $q \le n-1$,
\[
{E}_{p, q}^2=H_p(\fff, H_q(\SL_{n-1}, A))\simeq 
H_p(\fff, H_q(\SL, A))={E'}_{p, q}^2,
\]
which is induced by the pair
$(\id, \inc): (\fff, \SL_{n-1}) \arr (\fff, \SL)$. 

The spectral sequences give us a map of filtrations
\begin{gather*}
\begin{array}{ccccccc}
0=F_{-1} & \se & F_0 & \se \cdots \se & F_{n-1} 
&  \se & F_n=H_n(\GL_{n-1},A)\\
         &     & \downarrow&      & \downarrow  &      & \downarrow \\
0=F_{-1}'& \se & F_0'& \se \cdots \se & F_{n-1}'& \se  & F_n'=H_n(\GL, A),
\end{array}
\end{gather*}
such that ${E}_{i, n-i}^2\simeq F_i/F_{i-1}$, 
${E'}_{i, n-i}^2\simeq F_i'/F_{i-1}'$.
The commutative diagram
\[
\begin{array}{ccccccccc}
\!\!\! 0  & \!\!\!\larr &\!\!\!  F_{i-1}  &
\!\!\!
{\larr}
& \!\!\!  F_i &  \!\!\! \larr
&\!\!\! {E}_{i, n-i}^\infty
& \!\!\!\larr & \!\!\! 0 \\
&  &\!\!\!\Big\downarrow\vcenter{%
\rlap{$\scriptstyle{}$}}
&       & \!\!\!\Big\downarrow\vcenter{%
\rlap{$\scriptstyle{}$}}
&  &\!\!\!\Big\downarrow\vcenter{%
\rlap{$\scriptstyle{\wr}$}}
&      &  \\
\!\!\! 0  & \!\!\!\larr & \!\!\! F_{i-1}'
& \!\!\!\larr & \!\!\! F_i'
& \!\!\!  \larr   & \!\!\! {E'}_{i, n-i}^\infty    &\!\!\!
\larr & \!\!\! 0
\end{array}
\]
implies that
$\im (F_i \arr F_i') \cap F_{i-1}'= \im (F_{i-1} \arr F_{i-1}')$,
for $1 \le i \le n-1$.
From this we obtain
$\im (F_n \arr F_n') \cap F_{0}'= \im (F_{0} \arr F_{0}')$.
Since $H_n(\SL, A)\arr H_n(\GL, A)$ is injective,
$F_0'={E'}_{0, n}^\infty=H_n(\SL, A)$. In a similar way one obtains
$F_0={E}_{0, n}^\infty=H_n(\SL_{n-1}, A)_\fff/T$. It is easy to see that
the map $F_0 \arr F_0'$ is induced by the natural map
$\SL_{n-1} \arr \SL$ and so 
\[
T \se \ker (H_n(\SL_{n-1}, A)_\fff \arr H_n(\SL, A)).
\]
This completes the proof of (i). The proof of (ii) is 
analogue to (i).
\end{proof}

\begin{thm}\label{mir-exa}
{\rm (i)} If $B:=\nnn$, then
\[
H_{n}(\SL_{n-1}, B)_\fff
\arr H_{n}(\SL_{n}, B)_\fff
\arr K_n^M(F) \otimes B \arr 0
\]
is exact.
\par {\rm (ii)} If $F$ is  algebraically closed, then
\[
H_{n}(\SL_{n-1}) \arr H_{n}(\SL_{n}) \arr K_n^M(F) \arr 0
\]
is exact.
\end{thm}
\begin{proof}
 
(i) Consider the commutative diagram
\[
\begin{array}{ccccccc}
\!\!\!  H_n(\SL_{n-1}, B)_\fff  &
\!\!\!
{\larr}
& \!\!\!  H_n(\SL_{n}, B)_\fff &  \!\!\! \larr
&\!\!\! K_n^M(F)\otimes B
& \!\!\!\larr & \!\!\! 0 \\
\Big\downarrow\vcenter{%
\rlap{$\scriptstyle{}$}}
&       & 
\Big\downarrow\vcenter{%
\rlap{$\scriptstyle{}$}}
&  &
\Big\downarrow\vcenter{%
\rlap{$\scriptstyle{}$}}
&      &  \\
\!\!\! H_n(\GL_{n-1}, B)
& \!\!\!\larr & \!\!\! H_n(\GL_{n}, B)
& \!\!\!  \larr   & \!\!\! K_n^M(F)\otimes B   
&\!\!\! \larr & \!\!\! 0.
\end{array}
\]
Note that the second row in the above diagram is exact.
Suslin's map
$s_n: H_n(\GL_{n}) \arr K_n^M(F)$ has the property that the composition
$K_n^M(F) \arr K_n(F) \arr H_n(\GL_{n}) \arr K_n^M(F)$ coincides with
multiplication by $(-1)^{n-1}(n-1)!$. The Hurewicz map
$K_n(F) \arr H_n(\GL_{n})$ factors through 
$H_n(\SL_{n})_\fff \simeq H_n(\SL)$.  
Since $(n-1)!$ is invertible in $B=\nnn$,
the restriction map $s_n: H_n(\SL_{n}, B)_\fff \arr K_n^M(F)\otimes B$
is surjective. The exactness of the complex follows from
the exactness of the second row of the above diagram and Lemma \ref{inter}.

The proof of (ii) is analogue. In this case one should use
the fact that for an algebraically closed field $F$,
$K_n^M(F)$ is uniquely divisible \cite[1.2]{bass-tate1973}.
\end{proof}


Let $ K_n^M(F) \arr K_n(F)$ be the natural map from
the Milnor $K$-group to the Quillen $K$-group. Define
$K_n(F)^{\rm ind}:= \coker(K_n^M(F)\arr K_n(F))$.
This group is called the indecomposable part of $K_n(F)$.

\begin{cor}\label{k-ind}
Let $k$ be a field such that $(n-1)! \in k^\ast$. Assume
$H_n(\inc): H_n(\GL_{n-1}, k) \arr H_n(\GL_n, k)$ is injective. Then
\par {\rm (i)} $H_{n}(\SL_{n-1}, k)_\fff \overset{H_n(\inc)}{\larr}
H_{n}(\SL_{n}, k)_\fff$ is injective.
\par {\rm (ii)}
$K_n(F)^\ind \otimes k$ embeds in $H_{n}(\SL_{n-1}, k)_\fff$
and
\[
H_{n}(\SL_{n-1}, k)_\fff / K_n(F)^\ind \otimes k \simeq
H_{n}(\SL_{n}, k)_\fff / K_n(F)\otimes k.
\]
\par {\rm (iii)} In $H_n(\GL_{n}, k)$,
$K_n(F)\otimes k \cap H_n(\GL_{n-1}, k)$ coincides with
\[
K_n(F)\otimes k \cap  H_{n}(\SL_{n-1}, k)_\fff = K_n(F)^\ind\otimes k.
\]
\end{cor}
\begin{proof}
(i) This follows from the assumptions, \ref{mir-exa}, the commutativity of 
the diagram in the proof of \ref{mir-exa}  and the injectivity
of the natural map $\inc_\ast: H_{n}(\SL_{n-1}, k)_\fff \arr H_n(\GL_{n-1}, k)$.
For the proof of the injectivity of $\inc_\ast$
consider the map $\gamma_{n-1}: \fff \times \SL_{n-1} \arr \GL_{n-1}$
given by $(a, A) \mt aA$. By tracing the kernel and cokernel of
this map we obtain the exact sequence
\[
1 \arr \mu_{n-1, F} \arr \fff \times \SL_{n-1}
\overset{\gamma_{n-1}}{\arr} \GL_{n-1} \arr \fff/\fff^{n-1} \arr 1.
\]
From this we get two short exact sequences
\begin{gather*}
1 \arr \mu_{n-1, F} \arr \fff \times \SL_{n-1}
\arr \im(\gamma_{n-1}) \arr 1, \\
1 \arr \im(\gamma_{n-1})
\arr \GL_{n-1} \arr \fff/\fff^{n-1} \arr 1.
\end{gather*}
Writing the Lyndon-Hochschild-Serre spectral sequence of the above short
exact sequences (with coefficient in $k$)
and carrying out an easy analysis, one gets
\[
H_n(\fff \times \SL_{n-1}, k)_\fff \simeq H_n(\GL_{n-1}, k).
\]
Now by the K\"unneth theorem $H_n(\SL_{n-1}, k)_\fff$ embeds in 
$H_n(\GL_{n-1}, k)$, which  
%
is induced by the natural map 
$\inc: \SL_{n-1} \arr \GL_{n-1}$.\\
(ii) By (i) the complex
\[
0 \arr H_{n}(\SL_{n-1}, k)_\fff \arr H_{n}(\SL_{n}, k)_\fff
\arr K_n^M(F) \otimes k \arr 0
\]
is exact. By Suslin's construction of a map 
$K_n(F)\arr K_n^M(F)$ \cite[Section 4]{suslin1985},
we have a splitting exact sequence
\[
0 \arr K_n(F)^\ind\otimes k 
{\arr} K_n(F)\otimes k
{\arr} K_n^M(F)\otimes k \arr 0.
\]
The claims follow from applying the Snake lemma
to the commutative diagram
\[
\begin{array}{ccccccccc}
\!\!\! 0  & \!\!\!\larr &\!\!\!  K_n(F)^{\ind}\otimes k &
\!\!\!
{\larr}
& \!\!\!  K_n(F)\otimes k &  \!\!\! \larr
&\!\!\! K_n^M(F)\otimes k
& \!\!\!\larr & \!\!\! 0 \\
&  &\!\!\!\Big\downarrow\vcenter{%
\rlap{$\scriptstyle{{g}_n}$}}
&       & \!\!\!\Big\downarrow\vcenter{%
\rlap{$\scriptstyle{h_n}$}}
&  &\!\!\!\Big\downarrow  &      &  \\
\!\!\! 0  & \!\!\!\larr & \!\!\! H_{n}(\SL_{n-1}, k)_\fff
& \!\!\!\larr & \!\!\! H_{n}(\SL_{n}, k)_\fff
& \!\!\!  \larr   & \!\!\!K_n^M(F)\otimes k    &\!\!\!
\larr & \!\!\! 0,
\end{array}
\]
where $h_ n$ is the Hurewicz map
$K_n(F)\otimes k \arr H_n(\SL, k) \simeq H_{n}(\SL_{n}, k)_\fff$
which is injective \cite[Prop. 3, p. 507]{soule1985} and
${g}_n=h_n|_{K_n(F)^{\ind}}$.
\end{proof}

\section{Milnor $K$-groups}

We start this section with an easy lemma.
\begin{lem}\label{0-cycle}
Let $G$ be a group  and let
${g}_1, {g}_2,
h_1,\dots, h_n \in G$ such that each pair
commute. Let $C_G(\lan h_1,\dots, h_n \ran)$ be the subgroup
of $G$ consisting
of all elements of $G$ that commute with all $h_i$, $i=1, \dots, n$.
If ${\rm \bf{c}}({g}_1, {g}_2)=0$ in
$H_2(C_G(\lan h_1,\dots, h_n \ran))$, then
${\rm \bf{c}}({g}_1, {g}_2, h_1, \dots, h_n)=0$ in
$H_{n+2}(G)$.
\end{lem}
\begin{proof}
The homomorphism
$ C_G(\lan h_1,\dots, h_n \ran) \times \lan h_1,\dots, h_n \ran \arr G$
defined by $(g, h) \arr gh$ induces the
map
\[
H_2(C_G(\lan h_1,\dots, h_n \ran)) \otimes
H_n(\lan h_1,\dots, h_n \ran) \arr H_{n+2}(G).
\]
The claim follows from the fact that
${\rm \bf{c}}({g}_1, {g}_2, h_1, \dots, h_n)$ is the image of
${\rm \bf{c}}({g}_1, {g}_2)\otimes {\rm \bf{c}}(h_1, \dots, h_n)$
under this map.
\end{proof}

\begin{dfn}
Let $A_{i, n}:=\diag(a_i, \dots, a_i, a_i^{-(i-1)},I_{n-i}) \in \GL_n$.
We define
$[a_1, \dots, a_n]:={\bf c}(A_{1, n}, \dots , A_{n, n}) \in H_n(\GL_n)$.
\end{dfn}

\begin{prp}\label{kn}
{\rm (i)} The map $\nu_n: K_n^M(F) \arr H_n(\GL_n)$ defined by
\[
\{a_1, \dots, a_n\} \mt
[a_1, \dots, a_n]
\]
is a homomorphism.
\par {\rm (ii)} Let $s_n: H_n(\GL_n) \arr K_n^M(F)$
be the map defined by Suslin. Then 
$s_n \circ\nu_n$ coincides with the multiplication
by $(-1)^{(n-1)}(n-1)!$.
\end{prp}
\begin{proof}
(i)
It is well know that the Hurewicz map 
$h_2: K_2^M(F)\simeq K_2(F) \arr H_2(\GL)$
is defined by
$\{a, b\} \mt {\bf c}(\diag(a, 1, a^{-1}),\diag(b, b^{-1}, 1))$.
It is easy to see that in  $H_2(\GL)$,
$[a, b]= {\bf c}(\diag(a, 1, a^{-1}),\diag(b, b^{-1}, 1))$.
The stability isomorphism $H_2(\GL_2) \simeq H_2(\GL)$ implies that
the map $K_2^M(F)= K_2(F)\arr H_2(\GL_2)$, $\{a, b\} \mt [a, b]$
is well defined.
Now by lemma \ref{0-cycle},
$[a_1,1-a_1, a_3, \dots, a_n]=0$.
To complete the proof of (i) it is sufficient to prove that
\[
[a_1, \dots, a_{n-2},a_{n-1}, a_n]=-[a_1, \dots, a_{n-2}, a_n, a_{n-1}].
\]
This can be done in the following way;
\begin{gather*}
\begin{array}{l}
[a_1, \dots, a_{n-2}, a_{n-1}, a_n]\\
= {\bf c}(A_{1, n},\dots ,A_{n-2, n}, A_{n-1, n}, A_{n, n})\\
={\bf c}(A_{1, n}, \dots ,A_{n-2, n},\diag(a_{n-1}I_{n-2}, a_{n-1},
a_{n-1}^{-(n-1)}) , A_{n, n})\\
+{\bf c}(A_{1, n}, \dots , A_{n-2, n},
\diag(I_{n-2}, a_{n-1}^{-(n-1)},a_{n-1}^{(n-1)}) , A_{n, n})\\
 ={\bf c}(A_{1, n}, \dots,A_{n-2, n} ,
 \diag(a_{n-1}I_{n-2},a_{n-1},a_{n-1}^{-(n-1)}),\\
 \hspace{7.5 cm}
\diag(a_nI_{n-2},a_n^{-(n-2)}, 1))\\
+{\bf c}(A_{1, n}, \dots, A_{n-2, n} ,
\diag(a_{n-1}I_{n-2},a_{n-1},a_{n-1}^{-(n-1)}),\\
\hspace{7.5 cm}
\diag(I_{n-2},a_n^{(n-1)},a_n^{-(n-1)} ))\\
+{\bf c}(A_{1, n}, \dots,A_{n-2, n} ,
\diag(I_{n-2}, a_{n-1}^{-(n-1)},a_{n-1}^{(n-1)}) , A_{n, n})\\
=-[a_1, \dots, a_{n-2}, a_n,a_{n-1}]\\
+{\bf c}(A_{1, n}, \dots ,A_{n-2, n},
\diag(I_{n-2},a_{n-1},a_{n-1}^{-(n-1)}),\\
\hspace{7.5 cm}
\diag(I_{n-2},a_n^{(n-1)},a_n^{-(n-1)} ))\\
+{\bf c}(A_{1, n}, \dots ,A_{n-2, n}, \diag(a_{n-1}I_{n-2}, 1, 1),
\diag(I_{n-2},a_n^{(n-1)},a_n^{-(n-1)} ))\\
+{\bf c}(A_{1, n}, \dots ,A_{n-2, n},
\diag(I_{n-2}, a_{n-1}^{-(n-1)},a_{n-1}^{(n-1)}),
\diag(I_{n-2},a_n,a_n^{-(n-1)}))\\
+{\bf c}(A_{1, n}, \dots ,A_{n-2, n},
\diag(I_{n-2}, a_{n-1}^{-(n-1)},a_{n-1}^{(n-1)}),
\diag(a_nI_{n-2}, 1, 1))\\
=-[a_1, \dots, a_n,a_{n-1}].
\end{array}
\end{gather*}
(ii) Let $\tau_n$ be the composite map
$K_n^M(F) \arr K_n(F) \overset{h_n}{\arr} H_n(\GL_n)$. Then
$s_n\circ\tau_n$ coincides with the multiplication
by $(-1)^{(n-1)}(n-1)!$
\cite[section 4]{suslin1985}. It is well known that the composite map
$K_n^M(F)\overset{\tau_n}{\arr} H_n(\GL_n) 
\arr  H_n(\GL_n)/H_n(\GL_{n-1})$ is an isomorphism and it is defined by
$\{a_1, \dots, a_n\} \mt (a_1 \cup \dots \cup a_n) \mod H_n(\GL_{n-1})$, 
where
\[
a_1 \cup a_2 \cup \dots \cup a_n=
{\rm \bf{c}}(\diag(a_1, I_{n-1}),\diag(1, a_2, I_{n-2}),
\dots, \diag(I_{n-1},a_n))
\]
(see \cite[Remark 3.27]{nes-suslin1990}).
Also we know that $s_n$ factors as
\[
H_n(\GL_n) \arr  H_n(\GL_n)/H_n(\GL_{n-1}) \arr K_n^M(F).
\]
Our claim follows from the fact that modulo $H_n(\GL_{n-1})$
\[
[a_1, \dots , a_n]=(-1)^{n-1} (n-1)!
(a_1 \cup \dots \cup a_n).
\]
\end{proof}

\section{Homology of $\GL_n$}

Let $C_l(F^n)$ and $D_l(F^n)$ be the free abelian groups with
a basis consisting of $(\lan v_0\ran, \dots, \lan v_l\ran)$ and
$(\lan w_0\ran, \dots, \lan w_l\ran)$ respectively,
where every $\min\{l+1, n\}$
of $v_i \in F^n$ and every $\min\{l+1, 2\}$ of $w_i \in F^n$ are linearly
independent. By $\lan v_i\ran$ we mean the
line passing through  vectors $v_i$ and $0$.
Let $\partial_0: C_0(F^n) \arr
C_{-1}(F^n):=\z$, $\sum_i n_i(\lan v_i\ran) \mt \sum_i n_i$ and
$\partial_l=\sum_{i=0}^l(-1)^id_i: C_l(F^n) \arr C_{l-1}(F^n)$,
$l\ge 1$, where
\[
d_i((\lan v_0\ran, \dots, \lan v_l\ran))= (\lan
v_0 \ran, \dots,\widehat{\lan v_i \ran}, \dots, \lan v_l\ran).
\]
Define the differential $\tilde{\partial}_l=\sum_{i=0}^l(-1)^i
{\tilde{d}}_i: D_l(F^n)\arr D_{l-1}(F^n)$ similar to $\partial_l$.
Set $L_0=\z$, $M_0=\z$, $L_l=C_{l-1}(F^n)$ and 
$M_l=D_{l-1}(F^n)$, $l \ge 1$.
It is easy to see that the complexes
\begin{gather*}
\begin{array}{ll}
L_\ast \ : &
0 \leftarrow L_0 \leftarrow L_1  \leftarrow \cdots \leftarrow
L_l \leftarrow  \cdots \\
M_\ast: &
0 \leftarrow M_0 \leftarrow M_1 \leftarrow \cdots
\leftarrow M_l \leftarrow  \cdots
\end{array}
\end{gather*}
are exact. Take a $\GL_n$-resolution
$P_\ast \arr  \z$ of $\z$ with trivial $\GL_n$-action.
{}From the double complexes $L_\ast \otimes _{\GL_n} P_\ast$ and
$M_\ast \otimes _{\GL_n} P_\ast$ we obtain two first quadrant
spectral sequences converging to zero with
\begin{gather*}
\begin{array}{l}
E_{p, q}^1(n)= \begin{cases} H_q(\rr^p \times \GL_{n-p}) &
\text{if $0 \le p \le n$}\\
H_q(\GL_n, C_{p-1}(F^n) & \text{if $p\ge n+1$,}
\end{cases}\\
{\eee}_{p, q}^1(n)=
\begin{cases}
H_q(\fff^{p} \times \GL_{n-p})& \text{if $0 \le p \le 2$}\\
H_q(\GL_n, D_{p-1}(F^n))& \text{if $p \ge 3$.} \end{cases}
\end{array}
\end{gather*}
\\ For $1 \le p \le n$, and $q \ge 0$,  
$d_{p, q}^1(n)=
\sum_{i=1}^p(-1)^{i+1}H_q(\alpha_{i, p})$,
where 
\[
\begin{array}{l}
\alpha_{i, p}: \fff^p \times \GL_{n-p}\arr \fff^{p-1} \times
\GL_{n-p+1}, \\
\diag(a_1, \dots, a_p, A) \mt \diag(a_1, \dots,
\widehat{a_i}, \dots, a_p, \left(
\begin{array}{cc}
a_i & 0          \\
0   & A
\end{array}
\right)).
\end{array}
\]
In particular for $0 \le p \le n$, 
$d_{p, 0}^1(n)= \begin{cases} {\rm id}_\z & \text{if $p$ is odd}\\
0 & \text{if $p$ is even}\end{cases}$, so $E_{p, 0}^2(n)=0$ for $
p \le n-1$. It is also easy to see that $E_{n, 0}^2(n)=E_{n+1, 0}^2(n)=0$.
See the proof of \cite[Thm. 3.5]{mir1} for more details.

Let $k$ be a field and $C_i'(\fn):=C_i(\fn) \otimes k$.
Consider the commutative diagram of complexes
\begin{gather*}
\begin{array}{ccccccccccc}
0 & \leftarrow & k & \leftarrow &  C_0'(\fn)      & \leftarrow & C_1'(\fn)
& \leftarrow & C_2'(\fn) & \leftarrow & \cdots \\
  & & \Big\downarrow & & \Big\downarrow & & \Big\downarrow &
  & \Big\downarrow &  \\
0 & \leftarrow & 0 & \leftarrow & C_0'(\fn) & \leftarrow & C_1'(\fn)
& \leftarrow & C_2'(\fn) & \leftarrow & \cdots,
\end{array}
\end{gather*}
where the first vertical map is zero and the other vertical maps
are just identity maps.
This gives a map of the first quadrant spectral sequences
\begin{gather*}
E_{p, q}^1(n)\otimes k \arr \EE_{p, q}^1(n),
\end{gather*}
where  $\EE_{p, q}^1(n) \Rightarrow H_{p+q-1}(\GL_n, k)$
with $\EE^1$-terms
\[
\EE_{p, q}^1(n)=
 \begin{cases}
E_{p, q}^1(n)\otimes k &
\text{if $p \ge 1$}\\
0 & \text{if $p=0$} \end{cases}
\]
and differentials
$\DD_{p, q}^1(n)=
 \begin{cases}
d_{p, q}^1(n)\otimes \id_k &
\text{if $p \ge 2$}\\
0 & \text{if $p=1$} \end{cases}$.
It is not difficult to see that
$E_{p, q}^\infty \otimes k = \EE_{p, q}^\infty$
if $p\neq 1$, $q \le n$ and $p+q \le n+1$. Hence $\EE_{p, q}^\infty =0$
if $p\neq 1$, $q \le n$ and $p+q \le n+1$.

We look at the second spectral sequence in a different
way. The complex
\begin{gather*}
0  \leftarrow   C_0'(\fn)   \leftarrow  C_1'(\fn)
\leftarrow  \cdots \leftarrow  C_l'(\fn)  \leftarrow  \cdots
\end{gather*}
induces a  first quadrant spectral
sequence ${{\EE}'}_{p, q}^1(n) \Rightarrow H_{p+q}(\GL_n, k)$, where
${{\EE}'}_{p, q}^1(n)=\EE_{p+1, q}^1(n)$ and
${{\DD}'}_{p, q}^1(n)=\DD_{p+1, q}^1(n)$.
Thus ${{\EE}'}_{p, q}^\infty (n)=0$
if $p\ge 1$, $q \le n-1$ and $p+q \le n$.

\begin{prp}\label{exact-mir}
Let $n \ge 3$ and let $k$ be a field such that $(n-1)! \in k^\ast$.
Let the complex
\begin{equation*}\label{exact-k}
H_n(\fff^2\times \GL_{n-2}, k) \overset{\beta_2^{(n)}}{\larr}
H_n(\fff \times \GL_{n-1}, k) \overset{\beta_1^{(n)}}{\larr}
H_n(\GL_{n}, k) \arr 0
\end{equation*}
be exact, where $\beta_2^{(n)}=H_n(\alpha_{1,2})-H_n(\alpha_{2,2})$ and
$\beta_1^{(n)}=H_n(\inc)$. If the map
$H_m(\inc):H_m(\GL_{m-1}, k) \arr H_m(\GL_{m}, k)$ is
injective for $m=n-1, n-2$, then
$H_n(\inc): H_n(\GL_{n-1}, k) \arr H_n(\GL_n, k)$
is injective.
\end{prp}
\begin{proof}
The exactness of the above complex 
shows that the differentials
\[
{\DD'}_{r, n-r+1}^r(n): {\EE'}_{r, n-r+1}^r(n) \arr {\EE'}_{0, n}^r(n)
\]
are zero for $r \ge 2$. This proves that
${\EE'}_{0, n}^2(n)\simeq {\EE'}_{0, n}^\infty(n)$. To complete the proof
it is sufficient to prove that the group $H_n(\GL_{n-1}, k)$
is a summand of ${\EE'}_{0, n}^2(n)$. To prove this
it is sufficient to define a map 
\[
\varphi: H_n(\fff \times \GL_{n-1}, k) \arr
H_n(\GL_{n-1}, k)
\] 
such that
${\DD'}_{1, n}^1(H_n(\fff^2 \times \GL_{n-2}, k)) \se \ker(\varphi)$.
Consider the
decompositions $H_n(\fff \times \GL_{n}, k)= \bigoplus_{i=0}^n S_i$, where
$S_i= H_i(\fff, k) \otimes H_{n-i}(\GL_{n-1}, k)$.
For $2 \le i \le n$, the stability theorem gives the
isomorphisms $H_i(\fff, k) \otimes H_{n-i}(\GL_{n-2}, k) \simeq S_i$.
Define $\varphi: S_0 \arr H_n(\GL_{n-1}, k)$ the identity
map and for $2 \le i \le n$,
$\varphi: S_i \simeq H_i(\fff, k) \otimes H_{n-i}(\GL_{n-2}, k)
\arr H_n(\GL_{n-1}, k)$ the
shuffle product. To complete the definition of $\varphi$
 we must define it on $S_1$.
By a theorem of Suslin \cite[3.4]{suslin1985}
and the assumption, we have the decomposition
$H_{n-1}(\GL_{n-1}, k)\simeq H_{n-1}(\GL_{n-2}, k) 
\oplus K_{n-1}^M(F)\otimes k$.
So $S_1 \simeq H_1(\fff, k) \otimes H_{n-1}(\GL_{n-2}, k) \oplus
H_1(\fff, k) \otimes K_{n-1}^M(F)\otimes k$. Now define
$\varphi: H_1(\fff, k) \otimes H_{n-1}(\GL_{n-2}, k) 
\arr H_n(\GL_{n-1}, k)$
the shuffle product and
$\varphi: H_1(\fff, k) \otimes K_{n-1}^M(F) \arr H_n(\GL_{n-1}, k)$
the composite map
\begin{gather*}
\hspace{-2 cm}
H_1(\fff, k) \otimes K_{n-1}^M(F)\otimes k \overset{f}{\arr}
H_1(\fff, k) \otimes H_{n-1}(\GL_{n-1}, k)\\
\hspace{5 cm}
\overset{g}{\arr} H_n(\fff \times \GL_{n-1}, k) \overset{h}{\arr}
H_n(\GL_{n-1}, k),
\end{gather*}
where $f=\frac{1}{n-1}(\id \otimes \nu_{n-1})$,
$g$ is the shuffle product and
$h$ is induced by the map $\fff \times \GL_{n-1} \arr \GL_{n-1}$, 
$\diag(a, A) \mt aA$.
By the K\"unnuth theorem we have the decomposition
\begin{gather*}
\begin{array}{l}
T_0=H_{n}(\GL_{n-2}, k),\\
T_1=\bigoplus_{i=1}^n H_i(F_1^\ast, k)\otimes H_{n-i}(\GL_{n-2},k), \\
T_2=\bigoplus_{i=1}^n H_i(F_2^\ast , k)\otimes H_{n-i}(\GL_{n-2}, k), \\
T_3=H_1(F_1^\ast, k)\otimes H_1(F_2^\ast , k)\otimes
H_{n-2}(\GL_{n-2}, k),\\
T_4=\bigoplus_{ \substack {i+j \ge 3\\i, j \ne 0}}
H_i(F_1^\ast, k)\otimes H_j(F_2^\ast , k)\otimes H_{n-i-j}(\GL_{n-2}, k).
\end{array}
\end{gather*}
By lemma \ref{kn}, $T_3=T_3' \oplus T_3''$, where
\[
\begin{array}{l}
T_3'=H_1(F_1^\ast, k)\otimes H_1(F_2^\ast , k)
\otimes H_{n-2}(\GL_{n-3}, k),\\
T_3''=H_1(F_1^\ast, k)\otimes H_1(F_2^\ast , k)
\otimes K_{n-2}^M(F)\otimes k.
\end{array}
\]
It is not difficult to see that
${\DD'}_{1, n}^1(T_0 \oplus T_1\oplus T_2 \oplus T_3' \oplus T_4)
\se \ker(\varphi)$. Here one should use the stability theorem.
To prove ${\DD'}_{1, n}^1(T_3'') \se \ker(\varphi)$ we apply \ref{kn};
\begin{gather*}
\begin{array}{l}
{\DD'}_{1, n}^1\Big(a \otimes b \otimes \{c_1, \dots, c_{n-2}\}\Big)\\
=-\frac{(-1)^{n-3}}{(n-3)!}\bigg(
b \otimes {\bf c}(\diag(a, I_{n-2}),\diag(1, C_{1, n-2}), \dots ,
\diag(1, C_{n-2, n-2}))\\
\hspace{1.5 cm}
+a \otimes {\bf c}(\diag(b, I_{n-2}),\diag(1, C_{1, n-2}), \dots ,
\diag(1, C_{n-2, n-2}))\bigg)\\
=\frac{1}{(n-2)!}\bigg(
b \otimes [c_1, \dots, c_{n-2},a] + a \otimes [c_1, \dots, c_{n-2},b],\\
\hspace{1.5 cm}
-b \otimes{\bf c}(\diag(C_{1, n-2}, 1), \dots ,
\diag(C_{n-2, n-2}, 1), \diag(aI_{n-2}, 1))\\
\hspace{1.5 cm}
-a \otimes{\bf c}(\diag(C_{1, n-2}, 1), \dots ,
\diag(C_{n-2, n-2}, 1), \diag(bI_{n-2}, 1))\bigg).
\end{array}
\end{gather*}
Therefore
${\DD'}_{1, n}^1\Big(a \otimes b \otimes \{c_1, \dots, c_{n-2}\}\Big)=
(x_1, x_2) \in T_3' \oplus T_3''$, where
\begin{gather*}
\begin{array}{l}
x_1=-\frac{1}{(n-2)!}
\bigg(b \otimes{\bf c}(\diag(C_{1, n-2}), \dots ,
\diag(C_{n-2, n-2}), \diag(aI_{n-2}))\\
\hspace{2 cm}
+a \otimes{\bf c}(\diag(C_{1, n-2}), \dots ,
\diag(C_{n-2, n-2}), \diag(bI_{n-2}))\bigg),\\
x_2= (-1)^{n-2}
\bigg( b \otimes \{c_1, \dots, c_{n-2},a\} +
a \otimes \{c_1, \dots, c_{n-2},b\}\bigg).
\end{array}
\end{gather*}
We have $\phi(x_1)=-\frac{1}{(n-2)!}y$, where
\begin{gather*}
y=
{\bf c}(\diag(b,I_{n-2}), \diag(1, C_{1, n-2}), \dots ,
\diag(1, C_{n-2, n-2}), \diag(1, aI_{n-2}))\\
\hspace{0.4 cm}
+{\bf c}(\diag(a,I_{n-2}), \diag(1, C_{1, n-2}), \dots ,
\diag(1, C_{n-2, n-2}), \diag(1, bI_{n-2}))
\end{gather*}
and $\phi(x_2)=\frac{(-1)^{n-2}}{n-1}\frac{(-1)^{n-2}}{(n-2)!}z
=\frac{1}{(n-1)!}z$,
where
\begin{gather*}
\begin{array}{l}
z=\\
{\bf c}(\diag(bI_{n-1}), \diag(C_{1, n-2}, 1), \dots ,
\diag(C_{n-2, n-2}, 1), \diag(aI_{n-2},a^{-(n-2)}))\\
{\bf c}(\diag(aI_{n-1}), \diag(C_{1, n-2}, 1), \dots ,
\diag(C_{n-2, n-2}, 1), \diag(bI_{n-2},b^{-(n-2)}))\\
={\bf c}(\diag(bI_{n-1}), \diag(C_{1, n-2}, 1), \dots ,
\diag(C_{n-2, n-2}, 1), \diag(aI_{n-2},a))\\
+{\bf c}(\diag(bI_{n-1}), \diag(C_{1, n-2}, 1), \dots ,
\diag(C_{n-2, n-2}, 1), \diag(I_{n-2},a^{-(n-1)}))
\\
+{\bf c}(\diag(aI_{n-1}), \diag(C_{1, n-2}, 1), \dots ,
\diag(C_{n-2, n-2}, 1), \diag(bI_{n-2},b))\\
+{\bf c}(\diag(aI_{n-1}), \diag(C_{1, n-2}, 1), \dots,
\diag(C_{n-2, n-2}, 1), \diag(I_{n-2},b^{-(n-1)})).
\end{array}
\end{gather*}
Hence $\phi(x_2)=\frac{-1}{(n-2)!}z'$, where
\begin{gather*}
\begin{array}{l}
z'=\\
+{\bf c}(\diag(bI_{n-1}), \diag(C_{1, n-2}, 1), \dots ,
\diag(C_{n-2, n-2}, 1), \diag(I_{n-2},a))\\
+{\bf c}(\diag(aI_{n-1}), \diag(C_{1, n-2}, 1), \dots ,
\diag(C_{n-2, n-2}, 1), \diag(I_{n-2},b))\\
={\bf c}(\diag(bI_{n-2},1), \diag(C_{1, n-2}, 1), \dots ,
\diag(C_{n-2, n-2}, 1), \diag(I_{n-2},a))\\
+{\bf c}(\diag(I_{n-2},b), \diag(C_{1, n-2}, 1), \dots ,
\diag(C_{n-2, n-2}, 1), \diag(I_{n-2},a))\\
+{\bf c}(\diag(aI_{n-2},1), \diag(C_{1, n-2}, 1), \dots ,
\diag(C_{n-2, n-2}, 1), \diag(I_{n-2},b))\\
+{\bf c}(\diag(I_{n-2},a), \diag(C_{1, n-2}, 1), \dots ,
\diag(C_{n-2, n-2}, 1), \diag(I_{n-2},b))\\
=-y.
\end{array}
\end{gather*}
Therefore 
$\varphi(x_2)=\frac{-1}{(n-2)!}z'=-\frac{-1}{(n-2)!}y=-\varphi(x_1)$.
This completes the proof of the fact that
${\DD'}_{1, n}^1(H_n(\fff^2 \times \GL_{n-2}, k)) \se \ker(\varphi)$.
\end{proof}

Thus it is reasonable to conjecture

\begin{conj}\label{mir-conj}
Let $(n-1)! \in k^\ast$ and
let $n \ge 3$. Then
\begin{gather*}\label{exactness}
H_n(\fff^2\times \GL_{n-2}, k) \overset{\beta_2^{(n)}}{\larr}
H_n(\fff \times \GL_{n-1}, k)
\overset{\beta_1^{(n)}}{\larr}
H_n(\GL_n, k) \arr 0,
\end{gather*}
is exact.
\end{conj}

\begin{cor}
Let $(n-1)! \in k^\ast$.
If Conjecture $\ref{mir-conj}$ is true for all $n \ge 3$, then 
$H_n(\inc): H_n(\GL_{n-1}, k)\arr H_n(\GL_{n}, k)$ is 
injective for all $n$. In particular if $k=\q$, then 
Conjecture $\ref{mir-conj}$ implies Suslin's Injectivity Conjecture.
\end{cor}
\begin{proof}
This follows immediately from Proposition \ref{exact-mir}.
\end{proof}

\begin{rem}
(i) The surjectivity of $\beta_1^{(n)}$ is already  proven by Suslin
\cite{suslin1985}.
\par (ii) The conjecture is proven for n=3 in \cite[Prop. 2. 5]{mir2} and
we prove it in this note for $n=4$.
\par (iii) A similar result is not true
for $n=2$, that is
\begin{gather*}
H_2(\fff^2\times \GL_0) \overset{\beta_2^{(2)}}
{\larr}
H_2(\fff \times \GL_{1}) \overset{\beta_1^{(2)}}
{\larr} H_2(\GL_{2}) \arr 0
\end{gather*} 
is not exact. In fact
\[
\ker(\beta_1^{(2)}
)/\im(\beta_2^{(2)}
) \simeq \lan x \wedge  (x-1) -
x \otimes (x-1): x \in \fff \ran
\]
is a subset of $H_2(\fff) \oplus (\fff \otimes \fff)_\si$, where
$(\fff \otimes \fff)_\si=(\fff \otimes \fff)/
\lan a \otimes b + b \otimes a: a, b \in \fff \ran$.
To prove this let $Q(F)$ be the free abelian group with the
basis $\{[x]: x \in \fff-\{1\}\}$. Denote by $\ppp(F)$ the factor
group of $Q(F)$ by the subgroup generated by the elements of the form
$[x]-[y]+[y/x]-[(1-x^{-1})/(1-y^{-1})]+[(1-x)/(1-y)]$. The homomorphism
$\psi: Q(F) \arr \fff \otimes \fff,\ [x] \mt x \otimes (x-1)$
induces a homomorphism $\ppp(F) \arr (\fff \otimes \fff)_\si$,
\cite[1.1]{suslin1991}.
By \cite[2.2]{suslin1991}, $E_{4, 0}^2(2) \simeq \ppp(F)$.
It is not difficult to see that the
$E_{p, q}^2(2)$-terms have the following form
\begin{gather*}
\begin{array}{ccccccc}
\ast & \ast          &      &     &         &      &   \\
0    & E_{1, 2}^2(2) & \ast &     &         &      &   \\
0    & 0             &  0   &  0  & \ast    &      &   \\
0    & 0             &  0   &  0  & \ppp(F) & \ast & .
\end{array}
\end{gather*}
An easy calculation shows that
$E_{1, 2}^2(2) \se H_2(\fff) \oplus (\fff \otimes \fff)_\si$.
By \cite[2.4]{suslin1991}
$d_{4, 0}^3(2):E_{4, 0}^3(2) \arr E_{1, 2}^3(2)\simeq E_{1, 2}^2(2)$
is defined by $d_{4, 0}^3(2)([x])= x \wedge (x-1) - x \otimes (x-1)$.
Because the spectral sequence converges to zero 
we see that  $d_{4, 0}^3(2)$
is surjective and so $E_{1, 2}^2(2)$ is generated by elements
of the form $x \wedge (x-1) - x \otimes (x-1) \in
H_2(\fff) \oplus (\fff \otimes \fff)_\si$.
\end{rem}

\section{Homology of $\GL_4$}

\begin{lem}\label{e60}
$\eee_{p, 0}^2(4)$ is trivial for $0 \le p \le 6$.
\end{lem}

\begin{rem}
Lemma \ref{e60} gives a positive answer to a question asked by Dupont for
$n=3$ (see \cite[4.12]{sah1989}).
\end{rem}

\begin{lem}\label{e51}
$\eee_{p, 1}^2(4)$ is trivial for $0 \le p \le 5$.
\end{lem}

\begin{lem}\label{e42}
$\eee_{p, 2}^2(4)$ is trivial for $0 \le p \le 4$.
\end{lem}

\begin{lem}\label{e33}
$\eee_{p, 3}^2(4)$ is trivial for $0 \le p \le 3$.
\end{lem}

The proof of these lemmas is lengthy calculation. 
We prove them in the next section.

\begin{thm}\label{exact1}
The complex
\begin{gather*}
H_4(\fff^2\times \GL_2) \overset{\beta_2^{(4)}}{\arr}
H_4(\fff \times \GL_3) \overset{\beta_1^{(4)}}{\arr}
H_4(\GL_4) \arr 0
\end{gather*}
is exact.
\end{thm}
\begin{proof}
This follows from  \ref{e60}, \ref{e51}, \ref{e42}, \ref{e33} and
the fact that the spectral sequence converges to zero.
\end{proof}


\begin{thm}\label{g4-inj}
{\rm (i)} If  $\char(k) \neq 2, 3$, then
$
H_4(\GL_3, k) \overset{H_4(\inc)}{\larr} H_4(\GL_4, k)$ is injective.
\par {\rm (ii)} If $F$ be algebraically closed, then
$
H_4(\GL_3) \overset{H_4(\inc)}{\larr} H_4(\GL_4)$ is injective.
\end{thm}
\begin{proof}
(i) This follows from \ref{exact1},
\cite[Thm. 4.2]{mir2} and \ref{exact-mir}.\\
(ii) Since $F$ is algebraically closed, $H_2(F)$ and $K_2^M(F)$
are uniquely divisible, so
$\tors(F^\ast, H_2(F^\ast))=\tors(\fff, K_2^M(F))=0$.
Also from $H_2(\GL_3)\simeq H_2(\GL_2)=H_2(\fff) \oplus K_2^M(F)$ one
sees that $\tors(F^\ast, H_2(\GL_3))=0$. Therefore by the
K\"unneth theorem
\begin{gather*}
\begin{array}{l}
H_4(\fff \times \GL_3)\simeq \bigoplus_{i=0}^4 H_i(\fff)\otimes
H_{4-i}(\GL_3),\\
H_4(\fff^2 \times \GL_2)\simeq \bigoplus_{0 \le i+j \le 4}
H_i(F_1^\ast)\otimes H_j(F_2^\ast)\otimes H_{4-i-j}(\GL_2).
\end{array}
\end{gather*}
Now the proof of the claim is similar to the proof of 
\ref{exact-mir}, using
the fact that $H_2(\inc): H_2(\GL_1)\arr  H_2(\GL_2)$,
$H_3(\inc): H_3(\GL_2)\arr  H_3(\GL_3)$ are injective 
\cite[Thm. 4.2]{mir2}.
The only place where we need a modification is
the definition of the map $f$. First we define
\begin{gather*}
{\nu'}_3: K_3^M(F) \arr  H_3(\GL_3), \ \
\{a, b, c\} \mt [a^{1/3}, b, c].
\end{gather*}
This map is well defined as $K_3^M(F)$ is uniquely divisible. Set
$f:=\id \otimes {\nu'}_3: H_1(\fff) \otimes K_3^M(F)\arr
H_1(\fff) \otimes H_3(\GL_3)$. The rest of the proof  can be done
similar to the proof of \ref{exact-mir}. 
We leave the detail to the reader.
\end{proof}

\begin{exa}
$H_4(\inc): H_4(\GL_3(\R), \zzz) \arr H_4(\GL_4(\R), \zzz)$
is injective. It is well-known that
$K_3^M(\R)=\lan \{-1, -1, -1\} \ran \oplus V$, where
$\lan \{-1, -1, -1\} \ran$ is a group of order 2 
generated by $\{-1, -1, -1\}$
and $V$ is a uniquely divisible group. 
So the proof of our claim is similar to the proof of 
Theorem \ref{g4-inj}.
(We invert $2$ in the coefficient ring in order to
eliminate the 2-torsion elements that appear in the decomposition
of $H_4(\fff \times \GL_3)$ and $H_4(\fff^2 \times \GL_2)$. 
This might not be necessary.)
\end{exa}

\begin{cor}\label{mir-sus}
{\rm (i)}
If $\char(k) \neq 2, 3$, then
\[
0 \arr H_4(\GL_3, k) \arr H_4(\GL_4, k)
\arr K_4^M(F)\otimes k \arr 0.
\]
is split exact.
\par {\rm (ii)} If $F$ is algebraicly closed, then
\[
0 \arr H_4(\GL_3) \arr H_4(\GL_4)
\arr K_4^M(F) \arr 0.
\]
is split exact.
\end{cor}
\begin{proof}
The part (i) follows from \ref{g4-inj}, Suslin's exact sequence
(see the introduction).
For part (ii) we also need to know that $K_4^M(F)$ is uniquely divisible.
It is easy to see that  a splitting map can be defined by
\begin{gather*}
\{a, b, c, d \} \mt -\frac{1}{6} [a, b, c, d], \ \ \
\{a, b, c, d\} \mt [{a}^{-1/6}, b, c, d],
\end{gather*}
respectively.
\end{proof}

\begin{prp}\label{sl4-exact}
{\rm (i)}
If $\char(k) \neq 2, 3$, then
\[
0 \arr H_4(\SL_3, k)_\fff \arr H_4(\SL_4, k)_\fff
\arr K_4^M(F)\otimes k \arr 0.
\]
is split exact.
\par {\rm (ii)} If $F$ is algebraically closed, then 
\[
0 \arr H_4(\SL_3) \arr H_4(\SL_4)
\arr K_4^M(F) \arr 0.
\]
is split exact.
\end{prp}
\begin{proof}
The first part follows from \ref{mir-sus}, \ref{mir-exa} and \ref{k-ind}.
For the proof of (ii),  as in the proof of  \ref{k-ind}, 
it is sufficient to prove that the canonical
map $H_4(\SL_3) \arr H_4(\GL_3)$ is injective.
For this we look at the Lyndon-Hochschild-Serre spectral sequence
\[
E_{p, q}^2=H_p(\GL_3, H_q(\mu_{3, F}))
\Rightarrow H_{p+q}(\fff \times \SL_3)
\]
obtained from the extension
$1 \arr \mu_{3, F} \arr \fff \times \SL_3 \overset{\gamma_3}{\arr}
\GL_3 \arr 0$, where $\gamma_3(a, A):=a A$.
It is  well-known that $H_i(\fff)$, $K_i^M(F)$ and $K_i(F)$
are divisible \cite[Prop. 4.7]{dupont2001},
\cite[1.2]{bass-tate1973}, \cite{suslin1984}.
The Hurewicz map $K_i(F) \arr H_i(\SL)$ is isomorphism for $i=2,3$
(see \cite[5.2]{suslin1991} or
\cite[2. 5]{sah1989}). Since $B\GL^+ \sim B\fff \times B\SL^+$
\cite[5.3]{suslin1991}, 
one can see that $H_i(\GL_3)$ is divisible for $1 \le i \le 3$.  
Therefore the $E^2$-terms of the spectral sequence are of the 
following form
\begin{gather*}
\begin{array}{cccccc}
\z/3\z    & 0          &             &            &            &  \\
0         & 0          & 0           &            &            &  \\
\z/3\z    & 0          & \z/3\z      & \ast       &            &  \\
0         & 0          & 0           & 0          &            &  \\
\z/3\z    & 0          & \z/3\z      & 0          &  \ast      &  \\
\z        & H_1(\GL_3) & H_2(\GL_3)  & H_3(\GL_3) & H_4(\GL_3).&
\end{array}
\end{gather*}
An easy analysis of this spectral sequence shows that
\[
H_3(\gamma_3): H_3(\fff \times \SL_3) \arr H_3(\GL_3)
\]
is surjective with the kernel of order dividing 9.
By the K\"unneth theorem
\[
H_3(\fff \times \SL_3)\simeq H_3(\SL_3)\oplus \fff \otimes H_2(\SL_3)
\oplus H_3(\fff).
\]
Let $\xi \in \fff$ be the third root of unity, that is
$\xi^3=1$ and  $\xi \neq 1$. If we use
the bar resolution $C_\ast(G)$
to define the homology of a group $G$
(see \cite[p. 36]{brown1994}),
one can see that $\chi(\xi):=[\xi|\xi|\xi]+[\xi|\xi^2|\xi]
\in H_3(\fff)$ has order 3.
In fact in $C_\ast(G)_G$ ($G=\fff$)
\[
\partial_4([\xi|\xi|\xi|\xi]+[\xi|\xi^2|\xi|\xi]+[\xi|\xi|\xi|\xi^2]
+[\xi^2|\xi|\xi^2|\xi])=3\chi(\xi).
\]
In a similar way we can define
$\chi(\xi I_3) \in H_3(\SL_3)$ and $\chi(\xi I_3) \in H_3(\GL_3)$.
Now it is easy to see that
$H_3(\gamma_3)(\chi(\xi I_3), 0, 2\chi(\xi))=0$.
Therefore the kernel of $H_3(\gamma_3)$ is not trivial. Thus
$d_{4, 0}^2$ or $d_{4, 0}^4$ is trivial.
In either case this implies that
\[
H_4(\fff \times \SL_3)\arr H_4(\GL_3).
\]
is injective. Therefore $H_4(\SL_3)\arr H_4(\GL_3)$ is injective.

To give a splitting map one should note that in $H_4(\GL_4)$,
$[a, b, c, d]$ is equal to
\begin{gather*}
{\bf c}(\diag(a ,1, a^{-1}, 1), \diag(b, b^{-1}, 1, 1),
\diag(c ,1, c^{-1}, 1), \diag(d, d, d, d^{-3})).
\end{gather*}
We name this new version of $[a, b, c, d]$ by $[[a, b, c, d]]$.
Then $[[a, b, c, d]] \in H_4(\SL_4)$ and a 
splitting map can be defined by
\begin{gather*}
\{a, b, c, d\} \mt -\frac{1}{6} [[a, b, c, d]], \ \
\{a, b, c, d\} \mt [[a^{-1/6}, b, c, d]],
\end{gather*}
respectively.
\end{proof}

\begin{cor}
\par {\rm (i)}
If  $\char (k) \neq 2, 3$, then $K_4(F)^\ind \otimes k$ embeds
in the group $H_4(\SL_3, k)_\fff$.
\par {\rm (ii)} If $F$ is algebraically closed, then
$K_4(F)^\ind$ embeds in  $H_4(\SL_3)$.
\end{cor}
\begin{proof}
(i) This follows from \ref{g4-inj} and \ref{k-ind}. Part (ii) can be proven
in  a similar way using the facts that the Hurewicz map
$K_4(F) \arr H_4(\SL_4)$ is injective 
\cite[Thm. 7.23]{arlettaz2000} and
\[
0 \arr K_4^M(F) \arr  K_4(F) {\arr} K_4(F)^\ind \arr 0
\]
is split exact, for $K_4^M(F)$ \cite[1.2]{bass-tate1973} and
$K_4(F)$ \cite{suslin1984} are uniquely divisible.
\end{proof}

\section{Proof of lemmas \ref{e60}, \ref{e51}, \ref{e42}, \ref{e33}}
In this section we will assume that $n=4$. For simplicity
we will write $\eee_{p, q}^1$ in place of $\eee_{p, q}^1(4)$ and so on.
By $H_q(v)$ we mean $H_q(\stabe_{\GL_4}(v))$.
Here we look at $\fff \simeq H_1(\fff)$ as an abelian group
with multiplicative structure or as an abelian group with additive 
structure, depending on the context. 
This shall not cause any confusion. Note that under the isomorphism
$\fff \arr H_1(\fff)$ we have  $1 \mt 0$ and $a^{-1} \mt -a$. 

To prove
the lemmas \ref{e60}, \ref{e51}, \ref{e42} and \ref{e33}
we need to describe $\eee_{p, q}^1$  for $p=3, 4, 5$. Let
$w_i \in D_2(F^4)$,
$u_i \in D_3(F^4)$ and $v_i \in D_4(F^4)$, where
\begin{gather*}
w_1=(\lan e_1 \ran, \lan e_2 \ran, \lan e_3 \ran),\ \
w_2=(\lan
e_1 \ran, \lan e_2 \ran, \lan e_1+ e_2\ran),\\
\\
\!\!\!\!
u_1= (\lan e_1 \ran, \lan e_2 \ran,
\lan e_3 \ran, \lan e_4 \ran),
\ \ \ \ \ \ \ \
u_2= (\lan e_1 \ran, \lan e_2 \ran,
\lan e_3 \ran, \lan e_1+ e_2+e_3 \ran),\\
\hspace{-1 cm}
u_3= (\lan e_1 \ran, \lan e_2 \ran,
\lan e_3 \ran, \lan e_1+ e_2 \ran),\ \
u_4= (\lan e_1 \ran, \lan e_2 \ran,
\lan e_3 \ran, \lan e_2+e_3 \ran),\\
\hspace{-1 cm}
u_5= (\lan e_1 \ran, \lan e_2 \ran,
\lan e_3 \ran, \lan e_1+e_3 \ran), \ \
u_6= (\lan e_1 \ran, \lan e_2 \ran,
\lan e_1+e_2 \ran, \lan e_3\ran),\\
\hspace{-2.5 cm}
u_{7, a}= (\lan e_1 \ran, \lan e_2 \ran,
\lan e_1+e_2 \ran, \lan e_1+ae_2 \ran), \ \ \
a \in \fff-\{1\},\\
\\
\begin{array}{ll}
v_1= (\lan e_1 \ran, \lan e_2 \ran,
\lan e_3 \ran, \lan e_4 \ran, \lan e_1+ e_2+e_3 + e_4\ran) & \\
v_2^{i, j, k}= (\lan e_1 \ran, \lan e_2 \ran,
\lan e_3 \ran, \lan e_4 \ran, \lan e_i+ e_j+e_k\ran),&
1 \le i < j < k \le 4, \\
v_3^{i, j}= (\lan e_1 \ran, \lan e_2 \ran,
\lan e_3 \ran, \lan e_4\ran,\lan e_i+ e_j\ran),&
1 \le i < j \le 4, \\
v_4= (\lan e_1 \ran, \lan e_2 \ran,
\lan e_3 \ran,\lan e_1+ e_2+e_3\ran,\lan e_4\ran), &\\
v_{5, a, b}= (\lan e_1 \ran, \lan e_2 \ran,
\lan e_3 \ran, \lan e_1+ e_2+e_3 \ran, &
a, b \in \fff, a \neq 1 
\\
 &
\hspace{-3.1 cm}
\lan e_1+ ae_2+be_3\ran),\ {\rm or} \ b \neq 1 \\
v_{6, a}^{i, j}= (\lan e_1 \ran, \lan e_2 \ran,
\lan e_3 \ran, \lan e_1+ e_2+e_3\ran,\lan e_i+ ae_j\ran),& 
1 \le i < j \le 3, a \in \fff \\
v_7^{i, j}= (\lan e_1 \ran, \lan e_2 \ran,
\lan e_3 \ran, \lan e_i+ e_j\ran,\lan e_4\ran), &
1 \le i < j \le 3, \\
v_{8, a}^{i, j}= (\lan e_1 \ran, \lan e_2 \ran,
\lan e_3\ran,\lan e_i+ ae_j\ran,\lan e_1+ e_2 + e_3 \ran),& 
   1 \le i < j \le 3, a \in \fff 
\\
v_{9, a}^{i, j}= (\lan e_1 \ran, \lan e_2 \ran,
\lan e_3\ran,\lan e_i+ e_j\ran,\lan e_i+ ae_j \ran),&
   1 \le i < j \le 3, a \in \fff\\
v_{10}= (\lan e_1 \ran, \lan e_2 \ran,
\lan e_3\ran,\lan e_1+ e_2\ran,\lan e_1+ e_3 \ran),&\\
v_{11}= (\lan e_1 \ran, \lan e_2 \ran,
\lan e_3\ran,\lan e_1+ e_3\ran,\lan e_1+ e_2 \ran),&\\
v_{12}= (\lan e_1 \ran, \lan e_2 \ran,
\lan e_3\ran,\lan e_1+ e_2\ran,\lan e_2+ e_3 \ran),&\\
v_{13}= (\lan e_1 \ran, \lan e_2 \ran,
\lan e_3\ran,\lan e_2+ e_3\ran,\lan e_1+ e_2 \ran),&
\\
v_{14}= (\lan e_1 \ran, \lan e_2 \ran,
\lan e_3\ran,\lan e_1+ e_3\ran,\lan e_2+ e_3 \ran),&\\
v_{15}= (\lan e_1 \ran, \lan e_2 \ran,
\lan e_3\ran,\lan e_2+ e_3\ran,\lan e_1+ e_3 \ran),&\\
v_{16}= (\lan e_1 \ran, \lan e_2 \ran,
\lan e_1 +e_2 \ran,\lan e_3\ran,\lan e_4 \ran),&
\end{array}
\end{gather*}
\begin{gather*}
\begin{array}{ll}
v_{17, a}= (\lan e_1 \ran, \lan e_2 \ran,
\lan e_1 + e_2 \ran,\lan e_3 \ran,\lan e_1 +ae_2 + e_3 \ran),&
a \in \fff\\
v_{18, a}= (\lan e_1 \ran, \lan e_2 \ran,
\lan e_1 + e_2 \ran,\lan e_3 \ran,\lan e_1 +ae_2 \ran),&
a \in \fff-\{ 1\},\\
v_{19}= (\lan e_1 \ran, \lan e_2 \ran,
\lan e_1 + e_2 \ran,\lan e_3 \ran,\lan e_1 +e_3 \ran),&\\
v_{20}= (\lan e_1 \ran, \lan e_2 \ran,
\lan e_1 + e_2 \ran,\lan e_3 \ran,\lan e_2 +e_3 \ran),&\\
v_{21, a}= (\lan e_1 \ran, \lan e_2 \ran,
\lan e_1 + e_2 \ran,\lan e_1 + ae_2 \ran,\lan e_3 \ran),&
a \in \fff-\{ 1\},\\
v_{22, a, b}= (\lan e_1 \ran, \lan e_2 \ran,
\lan e_1 + e_2 \ran,\lan e_1 + ae_2 \ran, &
a, b \in \fff-\{ 1\}, a \neq b.\\
& \hspace{-2 cm}
 \lan e_1 + be_2 \ran),\\
\end{array}
\end{gather*}
By the Shapiro lemma
\begin{gather*}
\begin{array}{l}
\eee_{3, q}^1=H_q(w_1) \oplus H_q(w_2)\\
\eee_{4, q}^1= H_q(u_1)\oplus \cdots \oplus
H_q(u_{7, a})\\
\eee_{5, q}^1= H_q(v_1) \oplus \cdots \oplus
H_q(v_{22, a, b}).
\end{array}
\end{gather*}
So by the center killer lemma \cite[Thm. 1.9]{suslin1985}, 
we may assume that
\begin{gather*}
\eee_{3, q}^1=H_q(\fff^3 \times \GL_1) \oplus
H_q(\fff I_2 \times \GL_2)
\end{gather*}
and so on. Note that
\begin{gather*}
\begin{array}{ll}
{\tilde{d}}_{p, q}^1=d_{p, q}^1 \  {\rm for} \   p=1, 2,&
{\tilde{d}}_{3, q}^1|_{H_q(w_1)}=d_{3, q}^1, \\
{\tilde{d}}_{3, q}^1|_{H_q(w_2)}=H_q(\inc),&
{\tilde{d}}_{4, q}^1|_{H_q(u_1))}=(d_{4, q}^1, 0),\\
{\tilde{d}}_{4, q}^1|_{H_q(u_2)}=0,&
{\tilde{d}}_{4, q}^1|_{H_q(u_3)}=(-H_q(\inc), H_q(\inc)),\\
{\tilde{d}}_{4, q}^1|_{H_q(u_4)}=(-H_q(\inc), H_q(\alpha)), &
{\tilde{d}}_{4, q}^1|_{H_q(u_5)}=(-H_q(\beta), H_q(\gamma)),\\
{\tilde{d}}_{4, q}^1|_{H_q(u_6)}=(H_q(\inc), H_q(\inc)),&
{\tilde{d}}_{4, q}^1|_{H_q(u_{7, a})}=0,
\end{array}
\end{gather*}
where
\begin{gather*}
\alpha: (a, b, b, c) \mt (b, b, a, c),\\
\beta: (a, b, a, c) \mt (b, a, a, c),\\
\gamma: (a, b, a, c) \mt (a, a, b, c).
\end{gather*}
\begin{athm}[{\bf Lemma \ref{e60}}]
The group
$\eee_{p, 0}^2$ is trivial for $0 \le p \le 6$.
\end{athm}
\begin{proof}
The triviality of $\eee_{6, 0}^2$ is the most difficult
one, which we prove it here. The rest is much easier and we
leave it to the readers.\\
{\bf Triviality of $\eee_{6, 0}^2$}. The proof is in
four steps;\\
{ \bf Step 1}. The sequence
$0 \arr C_\ast(F^4)_{\GL_4} \arr D_\ast(F^4)_{\GL_4}
\arr Q_\ast(F^4){\GL_4} \arr 0$ is exact, where
$Q_\ast(F^4):= D_\ast(F^4)/C_\ast(F^4)$.\\
{\bf Step 2}. The group $H_5(Q_\ast(F^4)_{\GL_4})$ is trivial.\\
{\bf Step 3}. The map induced in homology by
$ C_\ast(F^4)_{\GL_4} \arr D_\ast(F^4)_{\GL_4}$ is
zero in degree $5$.\\
{\bf Step 4}. The group $\eee_{6, 0}^2$ is trivial.\\
{\bf Proof of step 1}. For
$i \ge -1$, $D_i(F^4)\simeq C_i(F^4)\oplus Q_i(F^4)$.
This decomposition is compatible with the action of
$\GL_4$, so we get an exact sequence of $\z[\GL_4]$-modules
\[
0 \arr C_i(F^4) \arr D_i(F^4)\arr Q_i(F^4) \arr 0
\]
which splits as a sequence of $\z[\GL_4]$-modules. One can easily
deduce the desired exact sequence from this. Note that this
exact sequence does not split as complexes.\\
{\bf Proof of step 2}. The complex $Q_\ast(F^4)$
induces a spectral sequence
\begin{gather*}
{\hat{E}}_{p, q}^1=
\begin{cases}
0& \text{if $0 \le p \le 2$}\\
H_q(\GL_4, Q_{p-1}(F^4))&
\text{if $p \ge 3$ } \end{cases}
\end{gather*}
which converges to zero. To prove the claim it is sufficient  to prove
that ${\hat{E}}_{6, 0}^2=0$ and to prove this it is sufficient to prove 
that ${\hat{E}}_{4, 1}^2={\hat{E}}_{3, 2}^2=0$. One can see that
\begin{gather*}
\begin{array}{l}
{\hat{E}}_{3,q}^1=H_q(w_2)\\
{\hat{E}}_{4,q}^1=H_q(u_2) \oplus \cdots \oplus
H_q(u_{7, a}),\\
{\hat{E}}_{5,q}^1=H_q(v_2^{i, j, k}) \oplus \cdots \oplus
H_q(v_{22, a, b}).
\end{array}
\end{gather*}
It is easy to see that
${\hat{d}}_{4,2}^1|_{H_2(u_6)}=H_2(\inc)$, so it
is surjective. Thus ${\hat{E}}_{3, 2}^2=0$.
Easy calculation shows that, $d \in \fff -\{ 1\}$ fixed,
\begin{equation}\label{11}
\hspace{-2.7 cm}
{\hat{d}}_{5,1}^1: H_1(v_{18, d}) \arr {\hat{E}}_{4,1}^1,
\ y \mt (0, y, 0, 0, y,-y),
\end{equation}
\begin{equation}\label{22}
\hspace{-3 cm}
{\hat{d}}_{5,1}^1: H_1(v_{3}^{1,3}) \arr {\hat{E}}_{4,1}^1,
\ \ y \mt (0,\ast, 0, y, 0, y),
\end{equation}
\begin{equation}\label{33}
\hspace{-3 cm}
{\hat{d}}_{5,1}^1: H_1(v_{9, a}^{1,2}) \arr {\hat{E}}_{4,1}^1,
\ \ y \mt (0, 0, 0, 0, 0, y),
\end{equation}
\begin{equation}\label{44}
\hspace{-0.4 cm}
{\hat{d}}_{5,1}^1: H_1(v_{3}^{2,3}) \arr {\hat{E}}_{4,1}^1,
\ \ (a, b, b, c) \mt (0,(b ,b ,c ,a),\ast, 0 , 0 , 0),
\end{equation}
\begin{equation}\label{55}
\hspace{-3 cm}
{\hat{d}}_{5,1}^1: H_1(v_4) \arr {\hat{E}}_{4,1}^1,
\ \ \ \ y \mt (y, 0 , 0 , 0 , 0 , 0),
\end{equation}
\begin{equation}\label{66}
{\hat{d}}_{5,1}^1: H_1(v_{3}^{3,4}) \arr {\hat{E}}_{4,1}^1, \ \ 
(0, a,0,0) \mt (0,(a ,0 ,0 ,-a), 0 , 0 , 0 , 0).
\end{equation}

Let $x=(x_2, x_3, \dots, x_{7, a}) \in \ker({\hat{d}}_{4,1}^1)$.
By  (\ref{11}), (\ref{22}), (\ref{33}), (\ref{44}) and (\ref{55})
we may assume $x_6=0$, $x_5=0$, $x_{7, a}=0$, $x_3=0$, $x_2=0$,
respectively.
If $x_4=(a, b, b, c)$, then ${\hat{d}}_{4,1}^1(x_4)=
(b, b, \mtx {a}{0}{0}{c})=0$.
Hence $b=0$ and $c=-a$. Now by (\ref{66}),
${\hat{d}}_{5,1}^1((0, a,0,0))=(0, 0, x_4, 0, 0, 0)=x$.
Therefore ${\hat{E}}_{4, 1}^2=0$.\\
{\bf Proof of step 3}. Here all the calculation will take place in
$
C_\ast(F^4)_{\GL_4}$ and
$
D_\ast(F^4)_{\GL_4}$. For simplicity
the image of $v \in C_\ast(F^4)$ in $C_\ast(F^4)_{\GL_4}$ is
denoted by $v$.
Consider the following commutative diagram
\begin{gather*}
\begin{array}{ccccc}
 C_6(F^4)_{\GL_4} & \arr &  C_5(F^4)_{\GL_4} &
\arr &  C_4(F^4)_{\GL_4} \\
\Big\downarrow & & \Big\downarrow & & \Big\downarrow  \\
D_6(F^4)_{\GL_4} & \arr & D_5(F^4)_{\GL_4} &
\arr &  D_4(F^4)_{\GL_4}  .
\end{array}
\end{gather*}
The generators of $C_5(F^4)_{\GL_4}$
are of the form $x_{a, b, c} $ with
\[
x_{a, b, c}=(\lan e_1 \ran, \lan e_2 \ran, \lan e_3 \ran, \lan e_4 \ran,
\lan e_1+e_2+e_3+e_4 \ran,\lan e_1+ ae_2+ be_3 +ce_4\ran),
\]
where $a, b, c \in \fff-\{1\}$ and $a \neq b$, $a \neq c$, $b \neq c$.
Since $C_4(F^4) \otimes_{\GL_4} \z=\z$,
$x_{a, b, c} \in \ker(\partial_5 )$ and
the elements of this form generate
$\ker(\partial_5 )$.
Hence  to prove this step it is sufficient to prove that
$x_{a, b, c} \in \im(\tilde{\partial}_6 )$. Let
\[
X_{a, b, c}=(\lan e_1 \ran, \lan e_2 \ran, \lan e_3 \ran, \lan e_4 \ran,
\lan e_1+e_2+e_3+e_4 \ran,\lan e_1+ ae_2+ be_3 +ce_4\ran,\lan e_1+ e_2\ran),
\]
where $a, b, c$ are as above. Then
\[
\tilde{\partial}_6 (X_{a, b, c} )=
v_{\frac{1-b}{1-c}, 1-b}-v_{\frac{a-b}{a-c}, a-b}+
u_{\frac{1-b}{a-b}}-u_{\frac{1-c}{a-c}}+u_{\frac{1}{a}}+ x_{a, b, c},
\]
where
\begin{gather*}
\begin{array}{l}
v_{g, h}=(\lan e_1 \ran, \lan e_2 \ran, \lan e_3 \ran, \lan e_4 \ran,
\lan e_1+e_2+e_3+e_4 \ran,\lan e_2+ ge_3 +he_4\ran),\\
u_{l}=(\lan e_1 \ran, \lan e_2 \ran, \lan e_3 \ran, \lan e_4 \ran,
\lan e_1+e_2+e_3+e_4 \ran,\lan e_1+ le_2\ran) \\
\hspace{0.6 cm}
-(\lan e_1 \ran, \lan e_2 \ran, \lan e_3 \ran, \lan e_4 \ran,
\lan e_1+e_2+e_3+e_4 \ran,\lan e_1+ e_2\ran),
\end{array}
\end{gather*}
$g, h, l \in \fff-\{1\}, g \neq h$.
So it is sufficient to prove that $v_{g, h}-v_{p, q}$ and
$u_l$ are in the image if $\tilde{\partial}_6 $.
Let
\begin{gather*}
\begin{array}{l}
U_l=(\lan e_1 \ran, \lan e_2 \ran, \lan e_3 \ran, \lan e_4 \ran,
\lan e_1+e_2+e_3+e_4 \ran,\lan e_1+ e_2\ran,\lan e_1+ le_2\ran),\\
V_l=(\lan e_1 \ran, \lan e_2 \ran, \lan e_3 \ran, \lan e_4 \ran,
\lan e_1+ e_2\ran,\lan e_1+ le_2\ran,\lan e_3+ e_4\ran ),
\end{array}
\end{gather*}
where $l \in \fff-\{1\}$. Then
$\tilde{\partial}_6 ((V_l-U_l) )=u_l$. Set
\[
T_{g, h}=(\lan e_1 \ran, \lan e_2 \ran, \lan e_3 \ran, \lan e_4 \ran,
\lan e_1+ e_2+e_3+ e_4\ran,\lan e_2+ge_3+ he_4\ran,\lan e_2+ e_3\ran),
\]
where $g, h \in \fff-\{1\}, g \neq h$. Then
\begin{gather*}
\begin{array}{l}
\tilde{\partial}_6 (T_{g, h}-T_{p, q})=
-s_{\frac{1}{1-h}}+s_{\frac{g}{g-h}}+s_{\frac{1}{1-q}}-s_{\frac{p}{p-q}}
-z_{\frac{1-h}{g-h}}+z_{\frac{1-q}{p-q}} \\
\hspace{2.9 cm}
+ y_{\frac{1}{g}}-y_{\frac{1}{p}}
+v_{g, h}-v_{p, q},
\end{array}
\end{gather*}
where
\begin{gather*}
\begin{array}{l}
s_a=(\lan e_1 \ran, \lan e_2 \ran, \lan e_3 \ran, \lan e_4 \ran,
\lan e_1+ e_2+e_3+ e_4\ran,\lan e_1+ae_3+ e_4\ran),\\
z_a=(\lan e_1 \ran, \lan e_2 \ran, \lan e_3 \ran, \lan e_4 \ran,
\lan e_1+ e_2+e_3+ e_4\ran,\lan e_2+ ae_3\ran)\\
\hspace{0.58 cm}
-(\lan e_1 \ran, \lan e_2 \ran, \lan e_3 \ran, \lan e_4 \ran,
\lan e_1+ e_2+e_3+ e_4\ran,\lan e_2+ e_3\ran),\\
y_a=(\lan e_1 \ran, \lan e_2 \ran, \lan e_3 \ran, \lan e_4 \ran,
\lan e_1+ e_2+e_3\ran,\lan e_1+ ae_2\ran)\\
\hspace{0.58 cm}
+(\lan e_1 \ran, \lan e_2 \ran, \lan e_3 \ran, \lan e_4 \ran,
\lan  e_2+e_3+ e_4\ran,\lan e_2+ ae_3\ran)\\
\hspace{0.58 cm}
-(\lan e_1 \ran, \lan e_2 \ran, \lan e_3 \ran, \lan e_4 \ran,
\lan  e_1+e_2+ e_3\ran,\lan e_1+ e_2\ran)\\
\hspace{0.58 cm}
-(\lan e_1 \ran, \lan e_2 \ran, \lan e_3 \ran, \lan e_4 \ran,
\lan  e_2+e_3+ e_4\ran,\lan e_2+ e_3\ran).
\end{array}
\end{gather*}
So to prove that $v_{g, h}-v_{p, q} \in
\im(\tilde{\partial}_6 )$  it is sufficient to prove that
$s_{a}-s_{ b}, a \neq b$, $z_a, y_a  \in \im(\tilde{\partial}_6 )$,
$a, b \in \fff-\{1\}$. Set
\begin{gather*}
\begin{array}{l}
Y_a=(\lan e_1 \ran, \lan e_2 \ran, \lan e_3 \ran, \lan e_4 \ran,
\lan e_1+ e_2+e_3\ran,\lan e_1+ae_2\ran, \lan e_1+e_2\ran),\\
Y_a'=(\lan e_1 \ran, \lan e_2 \ran, \lan e_3 \ran, \lan e_4 \ran,
\lan e_2+ e_3+e_4\ran,\lan e_2+ae_3\ran, \lan e_2+e_3\ran).
\end{array}
\end{gather*}
Then
$y_a=\tilde{\partial}_6(Y_a+Y_a')-2\tilde{\partial}_6(V_{\frac{1}{a}})$.
To prove that $z_a \in \im(\tilde{\partial}_6)$, set
\begin{gather*}
Z_a=(\lan e_1 \ran, \lan e_2 \ran, \lan e_3 \ran, \lan e_4 \ran,
\lan e_1+ e_2+e_3+ e_4\ran,\lan e_2+e_3\ran, \lan e_2+ae_3\ran).
\end{gather*}
By an easy calculation
$z_a=\tilde{\partial}_6(V_a)-\tilde{\partial}_6(Z_a)$.
If
\begin{gather*}
S_a=(\lan e_1 \ran, \lan e_2 \ran, \lan e_3 \ran, \lan e_4 \ran,
\lan e_1+ e_2+e_3+ e_4\ran,\lan e_1+ae_3+e_4\ran, \lan e_1+e_4\ran),
\end{gather*}
then $\tilde{\partial}_6(S_a-S_b)=R_{\frac{1}{1-a},\frac{1}{1-b}}+s_a-s_b$,
where, $a ,b \in \fff-\{1\}, a \neq b$ and
\begin{gather*}
\begin{array}{l}
R_{a, b}=(\lan e_1 \ran, \lan e_2 \ran, \lan e_3 \ran, \lan e_4 \ran,
\lan e_1+ e_2+ e_4\ran,\lan e_1+ ae_2+ e_4\ran)\\
\hspace{0.9 cm}
-(\lan e_1 \ran, \lan e_2 \ran, \lan e_3 \ran, \lan e_4 \ran,
\lan e_1+ e_2+ e_4\ran,\lan e_1+ be_2+ e_4\ran)\\
\hspace{0.9 cm}
-(\lan e_1 \ran, \lan e_2 \ran, \lan e_3 \ran, \lan e_4 \ran,
\lan e_2+ e_3+ e_4\ran,\lan e_2+ ae_3+ e_4\ran)\\
\hspace{0.9 cm}
+(\lan e_1 \ran, \lan e_2 \ran, \lan e_3 \ran, \lan e_4 \ran,
\lan e_2+ e_3+ e_4\ran,\lan e_2+ be_3+ e_4\ran).
\end{array}
\end{gather*}

Thus to prove that $s_a-s_b \in \im(\tilde{\partial}_6)$ it is sufficient
to prove that $R_{a, b} \in \im(\tilde{\partial}_6)$. Let
\begin{gather*}
\begin{array}{l}
Q_a=(\lan e_1 \ran, \lan e_2 \ran, \lan e_3 \ran, \lan e_4 \ran,
\lan e_1+ e_2+ e_4\ran,\lan e_1+ ae_2+ e_4\ran,\lan e_1+ e_4\ran )\\
\hspace{0.7 cm}
-(\lan e_1 \ran, \lan e_2 \ran, \lan e_3 \ran, \lan e_4 \ran,
\lan e_2+ e_3+ e_4\ran,\lan e_2+ ae_3+ e_4\ran,\lan e_1+ e_4\ran ),\\
N_{a, b}=(\lan e_1 \ran, \lan e_2 \ran, \lan e_3 \ran, \lan e_4 \ran,
\lan e_1+ e_4\ran,\lan e_1+ ae_4\ran )\\
\hspace{0.7 cm}
-(\lan e_1 \ran, \lan e_2 \ran, \lan e_3 \ran, \lan e_4 \ran,
\lan e_1+ e_4\ran,\lan e_1+ be_4\ran ),\\
P_a=(\lan e_1 \ran, \lan e_2 \ran, \lan e_3 \ran, \lan e_4 \ran,
\lan e_3+ e_4\ran,\lan ae_3+ e_4\ran ).
\end{array}
\end{gather*}
Then $\tilde{\partial}_6(Q_a-Q_b)=N_{\frac{1}{1-a},\frac{1}{1-b}}+
P_{\frac{1}{1-a}}
+P_{\frac{1}{1-b}}+R_{a, b}$. If
\begin{gather*}
\begin{array}{l}
O_{a, b}=(\lan e_1 \ran, \lan e_2 \ran, \lan e_3 \ran, \lan e_4 \ran,
\lan e_1+ e_4\ran,\lan e_1+ ae_4\ran ,\lan e_2+ e_3\ran)\\
\hspace{0.9 cm}
-(\lan e_1 \ran, \lan e_2 \ran, \lan e_3 \ran, \lan e_4 \ran,
\lan e_1+ e_4\ran,\lan e_1+ be_4\ran ,\lan e_2+ e_3\ran),\\
\ M_a=(\lan e_1 \ran, \lan e_2 \ran, \lan e_3 \ran, \lan e_4 \ran,
\lan e_3+ e_4\ran,\lan ae_3+ e_4\ran ,\lan e_1+ e_2\ran),
\end{array}
\end{gather*}
then $\tilde{\partial}_6(O_{a, b})=N_{a, b}$ and 
$\tilde{\partial}_6(M_{a})=P_{a}$.
This completes the proof of step 3.\\
{\bf Proof of Step 4}. writing the homological 
long exact sequence of the short
exact sequence obtained in the first step, we get the exact sequence
\[
H_5(C_\ast(F^4)_{\GL_4}) \arr H_5(D_\ast(F^4)_{\GL_4} )
\arr H_5(Q_\ast(F^4)_{\GL_4}).
\]
By steps 2 and 3,  $H_5(D_\ast(F^4)_{\GL_4})=0$,
but $\eee_{6, 0}^2=H_5(D_\ast(F^4)_{\GL_4})$. This completes the
proof of the triviality of $\eee_{6, 0}^2$.
\end{proof}

\begin{athm}[{\bf Lemma \ref{e51}}]
$\eee_{p, 1}^2(4)$ is trivial for $0 \le p \le 5$.
\end{athm}
\begin{proof}
The case $p=5$ is the  most difficult case and we show that
$\eee_{5, 1}^2(4)=0$. We leave the rest to the reader.
Let $x=(x_1, \dots, x_{22, a, b}) \in \ker({\tilde{d}}_{5,1}^1)$. Set
\begin{gather*}
V_{22, a, b}=(\lan e_1 \ran, \lan e_2 \ran, \lan e_3 \ran,
\lan e_1+ e_2\ran,\lan e_1+ ae_2\ran ,\lan e_1+ be_2\ran),
\end{gather*}
where $a, b \in \fff-\{1\}, a\neq b$. Then 
\begin{gather*}
{\tilde{d}}_{6,1}^1: H_1(V_{22, a, b}) \arr H_1(v_{9, a'}^{1,2})
\oplus H_1(v_{22, a, b}),\ \  y \mt (\ast, y).
\end{gather*}
So $x -{\tilde{d}}_{6,1}^1(x_{22, a, b})=(x_1', \dots, x_{21, a}',0)$.
So we may assume $x_{22, a, b}=0$. In a similar way we may assume
$x_{5, a, b}=x_{6, a}^{i, j}=
x_{8, a}^{i, j}=x_{9, a}^{i, j}=x_{10}=x_{11}=x_{12}=x_{13}=x_{14}=x_{15}=
x_{19}=x_{17, a}=x_{20}=x_7^{1,2}=x_{16}=0$.
Choose $a, b \in \fff-\{1\}$, $a \neq b$ and set
\begin{gather*}
\begin{array}{l}
{{}_1V}_7^{1,3}=(\lan e_1 \ran, \lan e_2 \ran,
\lan e_3\ran,\lan e_4\ran,\lan e_1 +e_4\ran,\lan e_2 +ae_4\ran),\\
{{}_2V}_7^{1,3}=(\lan e_1 \ran, \lan e_2 \ran,
\lan e_3\ran, \lan e_1 +e_3\ran,\lan e_4\ran,\lan e_1 +ae_3\ran),\\
{{}_2V}_7^{1,3}=(\lan e_1 \ran, \lan e_2 \ran,\lan e_1 +e_2\ran,
\lan e_1 +ae_2\ran,\lan e_3\ran,\lan e_1 +be_2\ran),\\
{{}_3V}_7^{1,3}=(\lan e_1 \ran, \lan e_2 \ran,\lan e_3\ran,
\lan e_2 +e_3\ran,\lan e_2 +ae_3\ran,\lan e_2 +be_3\ran).
\end{array}
\end{gather*}
Then
\begin{gather*}
\begin{array}{l}
{\tilde{d}}_{6,1}^1:
\bigoplus_{i=1}^{4} {{}_iV}_7^{1,3}
\arr H_1(v_3^{i, j})\oplus H_1(v_7^{1, 3})
\oplus H_1(v_{18, a}) \oplus H_1(v_{21, a}),\\
((a, b, c, b), -(b, c, b, a),-(b, b, c,a),(a, b, b, c)) \mt
(\ast,(b, c ,b ,a),\ast,\ast).
\end{array}
\end{gather*} 
So we may assume $x_7^{1, 3}=0$.
If
\begin{gather*}
{{}_1V}_2^{1,2,3}=(\lan e_1 \ran, \lan e_2 \ran,
\lan e_3\ran,\lan e_4\ran, \lan e_1 +e_2+e_3\ran),\lan e_1 +e_2\ran),\\
{{}_2V}_2^{1,2,3}=(\lan e_1 \ran, \lan e_2 \ran,
\lan e_3\ran,\lan e_4\ran, \lan e_1 +e_2+e_4\ran),\lan e_1 +e_2\ran),
\end{gather*}
then
\begin{gather*}
\begin{array}{l}
{\tilde{d}}_{6,1}^1: \bigoplus_{i=1}^{2} {{}_iV}_2^{1,2,3}
\arr H_1(v_3^{1,2})\oplus H_1(v_2^{1,2,3})\oplus H_1(v_2^{1,2,4}),\\
((a, a, a, b), (a, a, b, a)) \mt 
((2a, 2a, a + b, a +  b),-(a, a, a, b),(a, a, b, a)).
\end{array}
\end{gather*}
So we may assume $x_2^{1,2,3}=0$. In a similar way we may assume
$x_2^{2,3,4}=x_2^{1,3,4}=x_7^{1,2}=x_7^{1,3}=x_{16}=x_{18, a}=0$.
Therefore we reduce $x$ to
\[
x=(x_1, x_2^{1,2,4}, x_3^{i, j}, x_4, x_7^{2,3}, x_{21, a})
\in \ker({\tilde{d}}_{5,1}^1).
\]
Set
\[
\begin{array}{lll}
x_1=(a'', a'', a'', a''),& x_2^{1,2,4}=(s, t, s, s), &
x_3^{1,2}=(a_0, a_0, b_0, c_0), \\
x_3^{1,3}=(a_1, b_1, a_1, c_1),&
x_3^{1,4}=(a_2, b_2, c_2, a_2), & x_3^{2,3}=(a_3, b_3, b_3, c_3), \\
x_3^{2,3}=(a_4, b_4, c_4, b_4), & x_3^{3,4}=(a_5, b_5, c_5, c_5), &
x_4=s', s', s', t'), \\
x_{21, a}=(a', a', b', c'). & &
\end{array}
\]
Since ${\tilde{d}}_{5,1}^1(x)=0$, one gets
$x_7^{2,3}=0$, $a'=0$, $c'=-{b'}$, $s=-s'$, $t=-t'$, $a_1=b_3=b_4$.
Let $a \in \fff-\{1\}$ and fix $d \in \fff-\{1\}$. If
\begin{gather*}
\begin{array}{l}
V_{21, a}=(\lan e_1 \ran, \lan e_2 \ran,\lan e_1 + e_2 \ran,
\lan e_1 + ae_2 \ran,\lan e_3\ran,\lan e_4\ran),\\
V_{4}=(\lan e_1 \ran, \lan e_2 \ran,\lan e_3\ran,
\lan e_1 + e_2+e_3 \ran,\lan e_4\ran,\lan e_1 + e_2+d e_3 \ran), \\
V_2^{1,2,4}=(\lan e_1 \ran, \lan e_2 \ran,\lan e_3\ran,
\lan e_4\ran,\lan e_1 + e_2+e_3 \ran,\lan e_1 + e_2+d e_3 \ran),
\end{array}
\end{gather*}
then 
\begin{gather*}
\begin{array}{l}
{\tilde{d}}_{6,1}^1: H_1(V_{21, a}) \arr  H_1(v_{21, a}),\ \ 
(0, 0,-{b'}, 0)\mt ,(0, 0, b',-{b'}),\\
{\tilde{d}}_{6,1}^1: H_1(V_2^{1,2,4})\oplus H_1(V_{4}) \arr
H_1(v_2^{1,2,4})\oplus H_1(v_{4}),\\
\hspace{1 cm}
((s, s, s, t),(s, s, s, t)) \mt (-(s, s, t, s),(s, s, s, t)).
\end{array}
\end{gather*}
So we may assume $x_4= x_2^{1, 2, 4}= x_{21, a}=0$.
Let 
\[
V=(\lan e_1 \ran, \lan e_2 \ran,\lan e_3\ran,\lan e_4\ran,
\lan e_1 + e_2+ e_3 + e_4\ran,\lan e_2+ b e_4\ran).
\] 
Then
${\tilde{d}}_{6,1}^1: H_1(V) \arr H_1(v_1) \oplus H_1(v_3^{1,3})
\oplus H_1(v_3^{2,3})\oplus H_1(v_3^{2,4})$,
where $y=(a, a, a, a) \mt(y, y, y ,y)$.
In this way we may assume
$a_1=b_3=b_4=0$. In a similar way we may assume $a_0=c_5=0$.
If
\[
V_1=(\lan e_1 \ran, \lan e_2 \ran,\lan e_3\ran,\lan e_4\ran,
\lan e_1 + e_2+ e_3 + e_4\ran,\lan e_1+ b e_4\ran),
\]
then 
${\tilde{d}}_{6,1}^1: H_1(V_1) \arr H_1(v_1) \oplus H_1(v_3^{1,4})$,
$y=(a'', a'', a'', a'') \mt(-y, y)$. 
So we may assume $x_1=0$.
Thus 
$x=(x_3^{1,2}, x_3^{1,3},x_3^{1,4},x_3^{2,3},x_3^{2,4},x_3^{3,4})$,
where 
\begin{gather*}
\begin{array}{lll}
x_3^{1,2}=(0, 0, b_0, c_0),& x_3^{1,3}=(0, b_1, 0 ,c_1),&
x_3^{1,4}=(a_2, b_2, a_2, c_2),\\
x_3^{2,3}=(a_3, 0, 0, c_3),& x_3^{2,4}=(a_4, 0, c_4, 0),& 
x_3^{3,4}=(a_5, b_5, 0, 0).
\end{array}
\end{gather*}
From ${\tilde{d}}_{5,1}^1(x)=0$, we have the following relations
\begin{gather*}
\begin{array}{ll}
b_1+b_2+a_3-a_4+a_ 5=0, & 2b_1+c_2-c_4+b_5=0,\\
b_0+2c_1+a_2=0,         & c_0+c_1+c_3=0,\\
c_0-b_0-c_1+c_3=0,      & b_0-c_0- b_1+a_3=0,\\
-a_3+a_4-b_5-a_5=0,     & -c_3+c_4-a_5-b_5=0,\\
-b_1+b_2-c_2+c_4=0,     & -c_1+c_2-b_2+a_4=0.
\end{array}
\end{gather*}
Set 
\[
H_1(v_3^{i, j})=H_1(v_3^{1,2})\oplus H_1(v_3^{1,3})\oplus H_1(v_3^{1,4})
\oplus H_1(v_3^{2,3})\oplus H_1(v_3^{2,4})\oplus H_1(v_3^{3,4}).
\]
Let
\begin{gather*}
\begin{array}{l}
W_1=(\lan e_1 \ran, \lan e_2 \ran,\lan e_3\ran,
\lan e_1 + e_2+ e_3 \ran,\lan e_4\ran,\lan e_2+ e_4\ran),\\
W_2=(\lan e_1 \ran, \lan e_2 \ran,\lan e_3\ran,\lan e_4\ran,
\lan e_1 + e_2+ e_4 \ran,\lan e_2+ e_4\ran),\\
W_3=(\lan e_1 \ran, \lan e_2 \ran,\lan e_3\ran,
\lan e_4\ran,\lan e_1 + e_3+ e_4 \ran,\lan e_1+ e_4\ran),\\
W_4=(\lan e_1 \ran, \lan e_2 \ran,\lan e_3\ran,\lan e_4\ran,
\lan e_2 + e_3+ e_4 \ran,\lan e_2+ e_4\ran),\\
W_5=(\lan e_1 \ran, \lan e_2 \ran,\lan e_3\ran,
\lan e_4\ran,\lan e_1 + e_3+ e_4 \ran,\lan e_3+ e_4\ran),\\
W_6=(\lan e_1 \ran, \lan e_2 \ran,\lan e_3\ran,\lan e_4\ran,
\lan e_2 + e_3+ e_4 \ran,\lan e_3+ e_4\ran).
\end{array}
\end{gather*}
Then 
${\tilde{d}}_{6,1}^1:
H_1(W_1) \oplus H_1(W_2) \oplus H_1(V_4) \oplus H_1(V_2^{1,2,4})
\arr H_1(v_3^{i, j})$, where
\begin{gather*}
\hspace{-3.5 cm}
((0, 0, 0, c_3),(0, 0, c_3, 0),(0, 0, 0, c_3), (0, 0, 0, c_3)) \mt \\
\hspace{1 cm}
((0, 0, 0, c_3),-(0, c_3, 0, 0),(0, c_3, c_3, 0),-(0, 0, 0, c_3),
-(0, 0, c_3, 0), 0)
\end{gather*}
and $\bigoplus_{i=3}^6 H_1(W_i)\arr H_1(v_3^{i, j})$, where  
\begin{gather*}
\hspace{-3.5 cm}
((0, a_3, 0, 0),(a_3, 0, 0, 0),-(0, a_3, 0, 0),-(a_3, 0, 0, 0)) \mt \\
\hspace{1 cm}
(0, (0, a_3, 0, 0),(0, a_3, 0, 0),(a_3, 0, 0, 0),(a_3, 0, 0, 0),
-2(a_3, a_3, 0, 0)).
\end{gather*}
So we may assume $a_3=c_3=0$,that is
$x_3^{2,3}=0$. From the above equations we have $a_4=-c_4=a_5-b_5$.
Let
\begin{gather*}
\begin{array}{l}
{{}_1V}_3^{2,4}=(\lan e_1 \ran, \lan e_2 \ran,\lan e_3\ran,
\lan e_4\ran,\lan e_1 + e_3 \ran,\lan e_1+a e_3\ran),\\
{{}_1V}_3^{2,4}=(\lan e_1 \ran, \lan e_2 \ran,\lan e_3\ran,
\lan e_4\ran,\lan e_2 + e_4 \ran,\lan e_2+a e_4\ran),\\
{{}_1V}_3^{2,4}=(\lan e_1 \ran, \lan e_2 \ran,\lan e_3\ran,
\lan e_4\ran,\lan e_2 + e_3 \ran,\lan e_2+a e_3\ran),\\
{{}_1V}_3^{2,4}=(\lan e_1 \ran, \lan e_2 \ran,\lan e_3\ran,
\lan e_4\ran,\lan e_3 + e_4 \ran,\lan e_3+a e_4\ran),\\
\end{array}
\end{gather*}
Under the map ${\tilde{d}}_{6,1}^1:\bigoplus_{i=1}^4 H_1({}_iV_3^{2,4})
\arr H_1(v_3^{i, j})$ we have
\begin{gather*}
\hspace{-3.5 cm}
((0, b, 0, c), (c, 0, b, 0),(b, 0, 0, c),-(b, c, 0, 0)) \mt \\
\hspace{4 cm}
(0, 0, (0, b, c, 0), 0, (b-c, 0, c-b, 0),-(c, b, 0, 0)).
\end{gather*}
So we may assume $a_4=c_4=0$, that is $x_3^{2,4}=0$. If
\begin{gather*}
V_3^{3,4}=(\lan e_1 \ran, \lan e_2 \ran,\lan e_3\ran,
\lan e_4\ran,\lan e_1 + e_2 \ran,\lan e_3+e_4\ran),
\end{gather*}
then for
$H_1(V_3^{3,4}) \arr H_1(v_3^{i, j})$ we have
\[
(a, a, 0, 0) \mt (-(a, a, 0, 0), (a, a, 0, 0), 0, 0, 0, 0).
\]
So we may assume $a_5=b_5=0$, that is $x_3^{3,4}=0$. Therefore
$b_1=-c_1=b_2-c_2=b_0-c_0$. Let
\begin{gather*}
W_7=(\lan e_1 \ran, \lan e_2 \ran,\lan e_3\ran,
\lan e_4\ran,\lan e_1 + e_2+ e_3 \ran,\lan e_1+ e_3\ran),\\
W_8=(\lan e_1 \ran, \lan e_2 \ran,\lan e_3\ran,\lan e_4\ran,
\lan e_1 + e_3+ e_4 \ran,\lan e_1+ e_4\ran),\\
W_9=(\lan e_1 \ran, \lan e_2 \ran,\lan e_3\ran,
\lan e_4\ran,\lan e_1 + e_2+ e_3 \ran,\lan e_1+ e_2\ran),\\
W_{10}=(\lan e_1 \ran, \lan e_2 \ran,\lan e_3\ran,\lan e_4\ran,
\lan e_1 + e_3+ e_4 \ran,\lan e_1+ e_3\ran).
\end{gather*}
Then under the map
${\tilde{d}}_{6,1}^1:
\bigoplus_{i=7}^{10} H_1(W_i) \arr H_1(v_3^{i, j})$
we have
\begin{gather*}
\hspace{-3.5 cm}
z_1:=((0, 0, 0, a),-(0, a, 0, 0),-(0, 0, 0, a),(0, a, 0, 0)) \mt \\
\hspace{1.5 cm}
(-(0, 0, 2a, a),(0,-a, 0, a),(0,-a, a, 0), 0, (-a, 0, a, 0),(a, 2a, 0, 0)).
\end{gather*}
If $z_2=((0, 2a, 0, a),(a, 0, 2a, 0),(2a, 0, 0, a),-(2a, a, 0, 0))
\in \bigoplus_{i=1}^4 H_1({}_iV_3^{2,4})$, then
${\tilde{d}}_{6,1}^1 (z_1+z_2)=(-(0, 0, 2a, a),(0,-a, 0, a),
(0, a, 2a, 0), 0, 0, 0)$.
So we may assume $b_1=c_1=0$. Thus $a_2=c_2=0$ and $b_0=c_0=-a_2$.
If $V_3^{1,2}=(\lan e_1 \ran, \lan e_2 \ran,\lan e_3\ran,
\lan e_4\ran,\lan e_1 + e_4 \ran,\lan e_2+e_3\ran)$, then under
${\tilde{d}}_{6,1}^1: H_1(V_3^{1,2}) \arr H_1(v_3^{i, j})$ we have
$(b_0, 0, 0, b_0) \mt ((0, 0, b_0, b_0), 0, -(b_0, 0, 0, b_0), 0, 0)$.
This proves that $x \in \im({\tilde{d}}_{6,1}^1)$, and therefore
$\eee_{5,1}^2(4)$ is trivial.
\end{proof}

\begin{athm}[{\bf Lemma \ref{e42}}]
$\eee_{p, 2}^2$ is trivial for $0 \le p \le 4$.
\end{athm}
\begin{proof}
The difficult case is $p=4$. Let $x=(x_1, \dots, x_{7, a}) \in
\ker({\tilde{d}}_{4,2}^1)$. Consider the following maps
\begin{gather*}
\begin{array}{ll}
{\tilde{d}}_{4,2}^1: H_2(v_{18, a}) \arr H_2(u_3) \oplus H_2(u_6)
\oplus H_2(u_{7, a}),& y \mt (y, y, -y)\\
{\tilde{d}}_{4,2}^1
: H_2(v_{9, a}^{1,2}) \two H_2(u_{7, a}),& \\
{\tilde{d}}_{4,2}^1: H_2(v_3^{1,3}) \arr H_2(u_1) \oplus H_2(u_3)
\oplus H_2(u_5),& y \mt (\ast, \ast, y),\\
{\tilde{d}}_{4,2}^1: H_2(v_3^{2,3}) \arr H_2(u_1) \oplus H_2(u_3)
\oplus H_2(u_4), & y \mt (y, \ast, -y)\\
{\tilde{d}}_{4,2}^1: H_2(v_4) \arr H_2(u_1) \oplus H_2(u_2),& y \mt (0, y).
\end{array}
\end{gather*}
So we may assume $x_6=x_{7, a}=x_5=x_3=x_2=0$. Hence
$x=(x_1, x_4) \in H_2(u_1) \oplus H_2(u_4)$.
Let $H_2(u_1)=H_2(\fff^3 \times \GL_1)= \bigoplus_{i=1}^{10} T_i$ and
$H_2(u_4)=H_2(\fff I_2 \times \fff \times \GL_2)= 
\bigoplus_{i=11}^{16} T_i$,
where
\begin{gather*}
\begin{array}{llll} 
T_1=H_2(F_1^\ast),& T_2=H_2(F_2^\ast),& T_3=H_2(F_3^\ast),& 
\hspace{-0.3 cm}
T_4=H_2(F_4^\ast),\\
T_5=F_1^\ast \otimes F_2^\ast, & T_6=F_1^\ast \otimes F_3^\ast,&
T_7=F_1^\ast \otimes F_4^\ast,& 
\hspace{-0.3 cm}
T_8=F_2^\ast \otimes F_3^\ast,\\
T_9=F_2^\ast \otimes F_4^\ast,  &
T_{10}=F_3^\ast \otimes F_4^\ast,&
T_{11}=H_2(F^\ast),& 
\hspace{-0.3 cm}
T_{12}=H_2(F^\ast I_2), \\
T_{13}=H_2(\GL_1),&
T_{14}=F^\ast \otimes F^\ast I_2,& T_{15}=F^\ast \otimes \GL_1,&
\hspace{-0.3 cm}
T_{16}=F^\ast I_2 \otimes \GL_1. 
\end{array}
\end{gather*}
Set $x=(x_1, x_4)=(z_1, \dots, z_{16})$
with $z_i \in T_i$.
%
We look at the maps
\begin{gather*}
\hspace{-4 cm}
{\tilde{d}}_{4,2}^1: H_2(v_3^{1,2})=
H_2(\fff I_2 \times \fff \times \GL_1)
\arr  \bigoplus_{i=1}^{16} T_i,\\
(y_1, y_2, y_3, y_4, y_5, y_6) \mt
(y_1, y_1, y_2, y_3, 0, y_4, y_5, y_4, y_5, y_6, 0, 0, 0, 0, 0, 0).
\end{gather*}
and
\begin{gather*}
\hspace{-2 cm}
{\tilde{d}}_{4,2}^1: H_2(v_3^{3,4})=
H_2(\fff \times \fff \times \fff I_2)
 \arr \bigoplus_{i=1}^{16} T_i,\\
(y, y', 0) \mt (0, 0, 0, 0, y ,y' ,y' , 0 , 
0 , 0 , 0 , 0 , 0 , 0 , 0 , 0).
\end{gather*}
where
$H_2(v_3^{3,4})=H_2(\fff \times \fff \times \fff I_2)=
F_1^\ast \otimes F_2^\ast  \oplus  F_1^\ast\otimes \fff I_2 \oplus S$.
So we may assume $z_2=z_3=z_4=z_8=z_9=z_{10}=0$, $z_5=z_6=0$.
Now easy computation shows that
\begin{gather*}
\begin{array}{rrl}
-z_1 -z_{11}=0, & -z_{12}=0,& -z_{12}+z_{13}=0, \\
z_1 -z_{13}=0, & z_{14}=0. & 
\end{array}
\end{gather*}
Thus $z_1=z_{11}=z_{12}=z_{13}=z_{14}=0$. Now it is easy to see
that the rest of $z_i$'s should be zero. Thus  $\eee_{4,2}^1=0$.
\end{proof}

\begin{athm}[{\bf Lemma \ref{e33}}]
$\eee_{p, 3}^2$ is trivial for $0 \le p \le 3$.
\end{athm}
\begin{proof}
Let $\eee_{2, 3}^1=\bigoplus_{i=1}^{13} T_i$  and
$\eee_{3, 3}^1=\bigoplus_{i=1}^{26} S_i \oplus S'$
where
\begin{gather*}
\begin{array}{lll}
T_1=H_3(\GL_2), &  T_2=H_1(F_1^\ast)\otimes H_2(\GL_2),
& T_3=H_2(F_1^\ast)\otimes H_1(\GL_2),  \\
T_4=H_3(F_1^\ast), & T_5=H_1(F_2^\ast)\otimes H_2(\GL_2),
& T_6=H_2(F_2^\ast)\otimes H_1(\GL_2),\\
T_7=H_3(F_2^\ast),& T_8=H_1(F_1^\ast)\otimes H_2(F_2^\ast),
&T_9=H_2(F_1^\ast)\otimes H_1(F_2^\ast),
\end{array}
\\
\begin{array}{ll}
T_{10}=H_1(F_1^\ast)\otimes H_1(F_2^\ast)\otimes H_1(\GL_2),
& T_{11}=\tors(H_1(F_1^\ast), H_1(F_2^\ast)),\\
T_{12}=\tors(H_1(F_1^\ast),H_1(\GL_2) ),
& T_{13}=\tors(H_1(F_2^\ast),H_1(\GL_2)),
\end{array}
\end{gather*}
and
\begin{gather*}
\hspace{-0.6 cm}
\begin{array}{lll}
S_1=H_3(\GL_1),  & S_2=H_1(F_1^\ast)\otimes H_2(\GL_1),
&S_3=H_2(F_1^\ast)\otimes H_1(\GL_1),\\
S_4=H_3(F_1^\ast),&S_5=H_1(F_2^\ast)\otimes H_2(\GL_1),
& S_6=H_2(F_2^\ast)\otimes H_1(\GL_1),\\
S_7=H_3(F_2^\ast),& S_8=H_1(F_3^\ast)\otimes H_2(\GL_1),
& S_9=H_2(F_3^\ast)\otimes H_1(\GL_1),\\
S_{10}=H_3(F_3^\ast),& S_{11}=H_1(F_1^\ast)\otimes H_2(F_2^\ast),
& S_{12}=H_2(F_1^\ast)\otimes H_1(F_2^\ast),
\end{array}
\\
\begin{array}{ll}
S_{13}=H_1(F_1^\ast)\otimes H_2(F_3^\ast),
& S_{14}=H_2(F_1^\ast)\otimes H_1(F_3^\ast),\\
S_{15}=H_1(F_2^\ast)\otimes H_2(F_3^\ast),
& S_{16}=H_2(F_2^\ast)\otimes H_1(F_3^\ast),\\
S_{17}=H_1(F_1^\ast)\otimes H_1(F_2^\ast)\otimes H_1(\GL_1),
& 
S_{18}=H_1(F_1^\ast)\otimes H_1(F_3^\ast)\otimes H_1(\GL_1)
\\
S_{19}=H_1(F_2^\ast)\otimes H_1(F_3^\ast)\otimes H_1(\GL_1),
&
S_{20}=H_1(F_1^\ast)\otimes H_1(F_2^\ast)\otimes H_1(F_3^\ast),
\\
S_{21}=\tors(H_1(F_1^\ast), H_1(F_2^\ast)),
& S_{22}=\tors(H_1(F_1^\ast), H_1(F_3^\ast)),\\
S_{23}=\tors(H_1(F_1^\ast), H_1(\GL_1)),
&  S_{24}=\tors(H_1(F_2^\ast), H_1(F_3^\ast)), \\
S_{25}=\tors(H_1(F_2^\ast), H_1(\GL_1)),
& S_{26}=\tors(H_1(F_3^\ast), H_1(\GL_1)),\\
S'= H_3(F^\ast I_2 \times \GL_2).&
\end{array}
\end{gather*}
Note that by \cite[Prop. 3.1]{mir2} these decompositions are canonical.
Let $x=(x_1, \dots, x_{26}, x') \in \ker({\tilde{d}}_{3,3}^1)$.
Consider 
$H_3(v_1)=H_3(\fff^4)=\bigoplus_{i=1}^{26} S_i$.
We have ${\tilde{d}}_{4,3}^1: S_{17} \arr S_{17} \oplus S_{18}
\oplus S_{20}$ and ${\tilde{d}}_{4,3}^1: S_{19} \arr S_{19} \oplus S_{20}$
given by
\begin{gather*}
\begin{array}{l}
a \otimes b \otimes c \mt (-a \otimes b \otimes c,
b \otimes c \otimes a + a \otimes c \otimes b,-a \otimes b \otimes c)\\
a \otimes b \otimes c \mt (a \otimes b \otimes c,
-b \otimes c \otimes a- a \otimes c \otimes b-a \otimes b \otimes c).
\end{array}
\end{gather*}
repectivcely. So we may assume $x_{17}=x_{20}=0$. If
$x_{18}=a_{18}\otimes b_{18}\otimes c_{18}$,
$x_{19}=a_{19}\otimes b_{19}\otimes c_{19}$,
then $0={\tilde{d}}_{3,3}^1(x)=(z_1, \dots, z_{13}) \in
\bigoplus_{i=1}^{13} T_i$, where
\[
z_{10}=-a_{18}\otimes b_{18}\otimes \mtx {1} {0} {0} {c_{18}}
+ a_{19}\otimes b_{19}\otimes \mtx {1} {0} {0} {c_{19}}=0.
\]
Thus $x_{18}=x_{19}$. Let 
$H_3(u_6)=H_3(\fff I_2 \times \fff \times \GL_1)=
\bigoplus_{i=1}^{8} A_i \oplus A'$, where
\begin{gather*}
\begin{array}{lll}
A_1=H_3(F^\ast I_2), &
A_2=H_1(F^\ast)\otimes H_2(\GL_1),&
A_3=H_2(F^\ast)\otimes H_1(\GL_1),\\
A_4=H_3(F^\ast),&
A_5=H_1(F^\ast I_2)\otimes H_2(F^\ast),&
A_6=H_2(F^\ast I_2)\otimes H_1(F^\ast),
\end{array}
\\
\begin{array}{ll}
A_{7}=H_1(F^\ast I_2)\otimes H_1(F^\ast)\otimes H_1(\GL_1),&
A_{8}=\tors(H_1(F^\ast),H_1(\GL_1)).
\end{array}
\end{gather*}
Then ${\tilde{d}}_{4,3}^1: A_{7} \arr S_{18} \oplus S_{19}$,
$y \mt (y, y)$. So we may assume $x_{18}=x_{19}=0$.
Consider
${\tilde{d}}_{4,3}^1: S_{5} \arr S_{5} \oplus S_{9}
\oplus S_{13}\oplus S_{15}$ and
${\tilde{d}}_{4,3}^1: S_{6} \arr S_{6} \oplus S_{8}
\oplus S_{14}\oplus S_{16}$ given by
\begin{gather*}
\begin{array}{l}
a \otimes \sum [b | c] \mt (-a \otimes \sum [b | c],
-\sum[b | c] \otimes a ,a \otimes \sum[b | c], a \otimes \sum[b | c])\\
\sum[d | e] \otimes f \mt (-\sum[d | e] \otimes f,-f \otimes \sum[d | e],
\sum[d | e] \otimes f,\sum[d | e] \otimes f),
\end{array}
\end{gather*}
respectively. Thus we may assume that $x_5=x_6=0$. Now easy
calculation shows that $z_4=z_7=x_7+x_4'=0$, $x_4' \in S'$. Using the map
\begin{gather*}
{\tilde{d}}_{4,3}^1: A_{1} \arr S_{4} \oplus S_{7} \oplus S',
y \mt (-y,-y,\ast),
\end{gather*}
we may assume that $x_7=0$. Again consider the maps
\begin{gather*}
\begin{array}{ll}
A_{2} \oplus A_{3} \oplus A_{4}
\overset{{\tilde{d}}_{4,3}^1}{\arr} S_{8} \oplus
S_{9} \oplus S_{10}\oplus  S', &
(x_{8},x_{9}, x_{10}) \mt (x_{8},x_{9}, x_{10}, \ast), \\
A_{5} \oplus A_{6} \overset{{\tilde{d}}_{4,3}^1}{\arr}
S_{13} \oplus S_{14}  \oplus S_{15}\oplus S_{16}\oplus S',&
(x_{15},x_{16}) \mt (x_{15},x_{16}, x_{15},x_{16}, \ast).
\end{array}
\end{gather*}
So we may assume $x_{8}=x_{9}= x_{10}=x_{15}=x_{16}=0$.
Applying the maps
\begin{gather*}
\begin{array}{ll}
{\tilde{d}}_{4,3}^1: S_{25} \arr S_{22} \oplus S_{24}
\oplus S_{25}\oplus S_{26}, & y \mt (y, y,-y,-y)\\
{\tilde{d}}_{4,3}^1: S_{24} \arr S_{21} \oplus S_{24}
\oplus S_{25}, & y \mt (y,-y, 0)
\end{array}
\end{gather*}
one may assume $x_{25}=x_{24}=0$. Applying the map
${\tilde{d}}_{4,3}^1: A_{8} \arr T_{26} \oplus S'$, given by
$y \mt (y, \ast)$,
we may assume that $x_{26}=0$.
Let $H_3(u_4)=  \bigoplus_{i=1}^{8} B_i \oplus B'$, where
\begin{gather*}
\begin{array}{lll}
B_1=H_3(\GL_1), &
B_2=H_1(F^\ast)\otimes H_2(\GL_1),&
B_3=H_2(F^\ast)\otimes H_1(\GL_1),\\
B_4=H_3(F^\ast ),&
B_5=H_1(F^\ast)\otimes H_2(F^\ast I_2),&
B_6=H_2(F^\ast)\otimes H_1(F^\ast I_2),
\end{array}
\\
\begin{array}{ll}
B_{7}=\tors(H_1(F^\ast), H_1(F^\ast I_2)),&
B_{8}=\tors(H_1(F^\ast),H_1(\GL_1) ),
\end{array}
\end{gather*}
Using the map
${\tilde{d}}_{4,3}^1: B_{8} \arr S_{23} \oplus S'$,
$y \mt (y, \ast)$,
we may assume  $x_{23}=0$.
By easy calculation we have
$z_8=x_{11}-x_{13}=0$, $z_9=x_{12}-x_{14}=0$ and $z_{11}=x_{21}-x_{22}=0$.
Consider the following maps
\begin{gather*}
\hspace{-3 cm}
{\tilde{d}}_{4,3}^1: B_1 \oplus B_2 \oplus B_3 \oplus B_4
\arr S_{1} \oplus S_{2}
\oplus S_{3} \oplus S_{4}\oplus  S',\\
\hspace{6 cm}
(x_1, x_2,x_3,x_4) \mt (x_1, x_2,x_3,x_4, \ast),\\
\begin{array}{ll}
{\tilde{d}}_{4,3}^1: B_5 \oplus B_6 \arr S_{11} \oplus S_{12}
\oplus S_{7} \oplus S_{14}\oplus  S',&
(y, y') \mt (y, y',y, y', \ast),\\
{\tilde{d}}_{4,3}^1: B_{11} \arr S_{21} \oplus S_{22} ,&
y \mt (y, y).
\end{array}
\end{gather*}
So we may assume $x_1=x_2=x_3=x_4=x_{11}=x_{12} =x_{13} =x_{14}=
x_{21} =x_{22}=0$. This reduce $x$ to an element of the form
$x=x' \in S'$. But the map
${\tilde{d}}_{3,3}^1: H_3(w_2)=H_3(\fff I_2 \times \GL_2) \arr
\eee_{2,3}^1=H_3(\fff^2 \times \GL_2)$ is injective, thus
$x=0$. This completes the proof of the triviality of $\eee_{3,3}^2$.
\end{proof}


\bigskip

\address{{\footnotesize 
Fakult\"at f\"ur Mathematik,  Universit\"at Bielefeld, 
Postfach 100131,  D-33501 Bielefeld, Germany.\\
email:\ bmirzaii@math.uni-bielefeld.de
}}

\end{document}